\documentclass{amsart}
\usepackage{math}
\usepackage{psfrag,amsrefs}
\usepackage[all]{xy}

\renewcommand\germ{{\mathcal G}}

\newcommand\orbit{{\mathcal O}}
\newcommand\tile{{\mathcal T}}

\newcommand\tiling{{\mathfrak T}}

\newcommand\U{{\mathcal U}}
\newcommand\V{{\mathcal V}}
\newcommand\sqtile[4]{\xymatrix{{}\ar@{-}[d]_{#2}\ar@{-}[r]^{#3} & {}\ar@{-}[d]^{#4}\\ {}\ar@{-}[r]_{#1} & {}}}

\begin{document}
\title{Automata, Groups, Limit Spaces, and Tilings}
\author[L. Bartholdi]{Laurent Bartholdi}
\author[A. G. Henriques]{Andr\'e G. Henriques}
\author[V. V. Nekrashevych]{Volodymyr V. Nekrashevych}
\date{submitted 19 December 2004; last corrections 19 October 2005}
\begin{abstract}
  We explore the connections between automata, groups, limit spaces of
  self-similar actions, and tilings. In particular, we show how a
  group acting ``nicely'' on a tree gives rise to a self-covering of a
  topological groupoid, and how the group can be reconstructed from
  the groupoid and its covering. The connection is via finite-state
  automata. These define decomposition rules, or self-similar tilings,
  on leaves of the solenoid associated with the covering.
\end{abstract}
\maketitle

\tableofcontents

%%%%%%%%%%%%%%%%%%%%%%%%%%%%%%%%%%%%%%%%%%%%%%%%%%%%%%%%%%%%%%%%
\section{Introduction}
Groups generated by automata appeared in the early 60's, in particular
through the work of Aleshin~\cite{aleshin:burnside}; important
examples were studied later by Grigorchuk~\cite{grigorchuk:burnside}
and Gupta and Sidki~\cite{gupta-s:burnside}. A general theory of such
groups has only recently started to emerge.

The situation changed considerably when the last author introduced the
notion of ``limit space'' of a contracting group generated by
automata. This provided a bridge with dynamical systems, by
associating with such a group a space with covering map, and vice
versa.

This paper explores the degree to which these two notions are
equivalent. It turns out that the construction of a contracting group
from a limit space, and of a limit space from a contracting group, are
inverse to each other if the limit space is considered as a
topological orbispace (Theorem~\ref{thm:samegp}).

The orbispaces considered in this paper are quite complicated: for
example, they have uncountable fundamental group; the set of points
with non-trivial isotropy groups is dense, etc. Other constructions
would yield ``nicer'' limit orbispaces, e.g.\ such that the points
with non-trivial isotropy groups form a closed nowhere dense set; but
at the cost of making less transparent the similarity between
self-coverings of orbispaces and self-similar groups.

This paper also attempts to present in a uniform way automata, groups,
limit spaces, and the tilings that they carry. Many of the results are
not new, although wherever possible they are stated in greater generality.

The main connections, starting from an automaton ($\Pi$), or a
topological space with a branched covering ($X$), are as follows:
\[\xymatrix@!@+10ex{{\tiling} & %%BEFORE&: \ar@/^1pc/[d]^{\text{hull}}
  {\Pi}\ar@/^1pc/[r]^{\text{generate}}
  \ar@/_1pc/[d]_{\txt{\scriptsize inverse\\ \scriptsize paths}}
  \ar@{}[dr]|(0.35){\displaystyle\dagger}
  \ar@{}[dr]|(0.75){\displaystyle\ddagger} &
  {\Gamma}\ar@/^1pc/[d]^{\text{complete}}\ar@/^1pc/[l]^{\text{nucleus}} \\
  {S}\ar@/^1pc/[u]^{\text{leaves}}\ar@/_1pc/[r]_{S/[x,y]\sim[x,y']} & 
  {L}\ar@/_1pc/[r]_{\text{Galois}}\ar[ur]_{\pi_1(L)/\text{kernel}}
  %%ADD_DISKS: \ar@/_1pc/[d]_{\text{add disks}}
  \ar@/_1pc/[l]_{\varprojlim} &
  {\overline\Gamma}\\
  & {X}\ar@/_1pc/[u]_{\txt{\scriptsize remove rami-\\ \scriptsize
  fication orbit}}}
\]
The connecting constructions are a limit space ($L$), a Solenoid ($S$)
and a discrete group ($\Gamma$). The commutativity of the diagram at
$\dagger$ is proven in Theorem~\ref{thm:samegp}, and that at
$\ddagger$ is proven in Theorem~\ref{thm:sameclosure}.

%%%%%%%%%%%%%%%%%%%%%%%%%%%%%%%%%%%%%%%%%%%%%%%%%%%%%%%%%%%%%%%%
\section{Definitions}
We introduce in this section the main definitions used in the
text. Most of them already appeared in~\cite{bartholdi-g-n:fractal}
with examples and illustrations.

Our convention for $\N$ is that it does not contain $0$.

\subsection{Automata}\label{subs:automata} An \emph{automaton} $\Pi$
is a pair $A,Q$ of sets called respectively \emph{alphabet} and
\emph{states}, with maps $\sigma:A\times Q\to A$ and $\tau:A\times
Q\to Q$ called respectively the \emph{output} and \emph{transition}
functions.

A \emph{graph} is a pair $V,E$ of sets called respectively
\emph{vertices} and \emph{edges}, with maps $s,t:E\to V$ called
respectively the \emph{source} and \emph{target} maps.\footnote{By
  ``graph'', we shall therefore always mean ``oriented graph'', with
  both loops and multiple edges allowed.}

The \emph{graph} $\gf(\Pi)$ of an automaton $\Pi$ is the graph with
vertex set $Q$, containing for all $q\in Q$ and all $a\in A$ an edge
from $q$ to $\tau(a,q)$, labeled $a/\sigma(a,q)$.

The \emph{dual} of $\Pi$ is the automaton $\Pi^*$ with alphabet $Q$,
states $A$, output $\sigma^*(q,a)=\tau(a,q)$ and transition
$\tau^*(q,a)=\sigma(a,q)$.

Let $\Pi$ and $\Pi'$ be two automata on the same alphabet $A$. The
\emph{product} of $\Pi$ and $\Pi'$ is the automaton $\Pi''=\Pi*\Pi'$ with
alphabet $A$, states $Q\times Q'$, output
$\sigma''(a,(q,q'))=\sigma'(\sigma(a,q),q')$ and transition
$\tau''(a,(q,q'))=(\tau(a,q),\tau'(\sigma(a,q),q'))$.
The iterated product $\Pi*\dots*\Pi$ is written $\Pi^n$.

We denote by $A^*=\cup_{n\ge0}A^n$ the free monoid on $A$. It is
conveniently represented as an $\#A$-regular rooted tree, with root the
empty word, and with an edge from $w$ to $wa$ for all $w\in A^*,a\in A$.

The output and transition functions are naturally extended to
functions $\sigma:A^*\times Q^*\to A^*$ and $\tau:A^*\times Q^*\to
Q^*$, by
\begin{align*}
  \sigma(a_1a_2\dots a_n,q)&=\sigma(a_1,q)\sigma(a_2\dots a_n,\tau(a_1,q)),\\
  \tau(a_1a_2\dots a_n,q)&=\tau(a_2\dots a_n,\tau(a_1,q)),\\
  \sigma(a_1\dots a_n,q_1q_2\dots q_m)&=\sigma(\sigma(a_1\dots
  a_n,q_1),q_2\dots q_m),\\
  \tau(a_1\dots a_n,q_1q_2\dots q_m)&=\tau(a_1\dots
  a_n,q_1)\tau(\sigma(a_1\dots a_n,q_1),q_2\dots q_m).
\end{align*}

The extended input and output functions can be visualized in the
following way: consider for all $a\in A,q\in Q$ a small square with
lower side labeled $a$, left side $q$, top side $\sigma(a,q)$ and
right side $\tau(a,q)$.  Consider a large rectangle $R$ with bottom
label $a_1\dots a_n$ and left label $q_1\dots q_m$, and tile it by the
above small squares. Then $\sigma(a_1\dots a_n,q_1q_2\dots q_m)$
and $\tau(a_1\dots a_n,q_1q_2\dots q_m)$ are the top and right labels
of $R$.

\begin{center}
  \psfrag{a}{$a$}
  \psfrag{q}{$q$}
  \psfrag{s}{$\sigma(a,q)$}
  \psfrag{t}{$\tau(a,q)$}
  \psfrag{A}{\hbox to 4in{$a_1$\dotfill $a_n$}}
  \psfrag{Q}[Bl][Bl][1][90]{\hbox to 2in{$q_1$\dotfill $q_m$}}
  \psfrag{S}{\hbox to 4in{\dotfill$\sigma(a_1\dots a_n,q_1\dots q_m)$\dotfill}}
  \psfrag{T}[Bl][Bl][1][90]{\hbox to 2in{\dotfill$\tau(a_1\dots
  a_n,q_1\dots q_m)$\dotfill}}
  \psfrag{R}{\LARGE $R$}
  \includegraphics{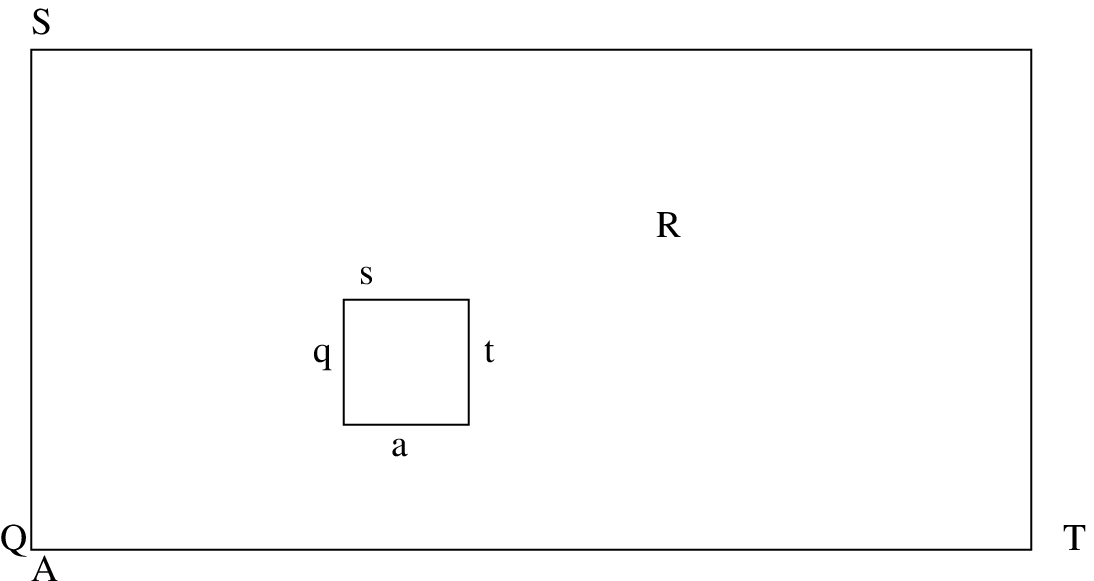}
\end{center}

Let $\Pi$ and $\Pi'$ be two automata on the same alphabet $A$. An
\emph{automaton homomorphism} $f:\Pi\to\Pi'$ is a map $f:Q\to Q'$ such
that $\sigma'(a,f(q))=\sigma(a,q)$ and $\tau'(a,f(q))=f(\tau(a,q))$
for all $a\in A,q\in Q$. Note that there is a label-preserving graph
homomorphism $\gf(\Pi)\to\gf(\Pi')$ if and only if there is an
automaton homomorphism $\Pi\to\Pi'$.  We call $\Pi'$ a
\emph{subautomaton} of $\Pi$ if $f$ is injective and $\Pi'$ is a
\emph{quotient} of $\Pi$ if $f$ is surjective.

\subsection{Actions} An automaton $\Pi$, with a given state $q$,
induces a transformation $\Pi_q$ on the tree $A^*$, via $\sigma$;
explicitly,
\[\Pi_q(a_1a_2\dots a_n)=\sigma(a_1,q)\,\Pi_{\tau(a_1,q)}(a_2\dots a_n).\]
Note that $\Pi'_{q'}\Pi_q=(\Pi*\Pi')_{(q,q')}$.

Every orbit of this action lies in $A^n$ for some
$n\ge0$.\footnote{Equivalently, the action is \emph{radial}.}

The transformations $\Pi_q$ are invertible if $\sigma(-,q)$ is a
permutation for all $q\in Q$.  In that case $\Pi$ is called
\emph{invertible}, and $\langle \Pi\rangle$, the \emph{group of the
  automaton $\Pi$}, is defined as the group generated by the $\Pi_q$
for all $q\in Q$.

Let $\Gamma$ be a group acting on a set $X$, with generating set $S$.
The \emph{Schreier graph} of $X$ is the graph with vertex set $X$, and
for all $x\in X$ and $s\in S$ an edge from $x$ to $sx$, labeled $s$.

\begin{lem}
  The Schreier graph of $\langle\Pi\rangle$ on $A^n$ is
  $\gf((\Pi^*)^n)$.
\end{lem}
\begin{proof}
  Direct verification. Indeed the vertices of $\gf(\Pi^*)$ naturally
  identify with $A$, and its edges identify with $A\times Q$.
  Similarly, the vertices in $\gf((\Pi^*)^n)$ identify with $A^n$, and
  its edges identify with $A^n\times Q$.
\end{proof}

\begin{defn}
  A group $\Gamma$ acting on a rooted tree $A^*$ is \emph{spherically
    transitive} if for every $n$ the action is transitive on $A^n$; in
  other words, $\Gamma$ acts transitively on the spheres around the
  root vertex.
\end{defn}

\begin{defn}
  A group $\Gamma$ is \emph{recurrent} if for any vertex $w\in A^*$
  the natural map $\stab_\Gamma(w)\to\aut(A^*)$, mapping $g\in\Gamma$
  to its restriction on $wA^*$ and identifying $wA^*$ with $A^*$, has
  image $\Gamma$.  Note that if $\Gamma$ is generated by an automaton,
  then this image is by necessity a subgroup of $\Gamma$; in that
  case, the statement needs to be checked only for all $w\in A$.
\end{defn}

An automaton $\Pi$ is \emph{recurrent} or \emph{spherically transitive}
if the associated group $\langle\Pi\rangle$ enjoys the respective
property.

\begin{lem}
  If $\Gamma$ is recurrent and acts transitively on $A$, then it is
  spherically transitive.
\end{lem}

\subsection{Contraction}
Let $\Pi$ be an automaton. It is \emph{contracting} if there is a
finite set $N\subset\langle\Pi\rangle$ such that for any
$g\in\langle\Pi\rangle$ there is $n\in\N$ with $\tau(A^n,g)\subset
N$. The minimal such set $N$ is called the \emph{nucleus} of
$\langle\Pi\rangle$. It may be defined as
\[N =
\bigcup_{g\in\langle\Pi\rangle}\bigcap_{n_0\ge0}\bigcup_{n\ge n_0}\tau(A^n,g).\]

\begin{defn}
  The automaton is \emph{nuclear} if it is contracting and its nucleus
  is equal to $Q$.
\end{defn}

We remark that it is usually difficult to determine whether an
automaton is contracting. On the other hand, it takes polynomial time
to check whether it is nuclear.

If $\Pi$ is a finite, contracting automaton, and is generated by its
nucleus $N$, then we may replace $\Pi$'s set of states with $N$, and
obtain a nuclear automaton generating the same group.

\begin{lem}\label{lem:tauepi}
  If $\Pi$ be a nuclear automaton, then $\tau:A\times Q\to Q$ is onto.
\end{lem}
\begin{proof}
  Assume to the contrary that $q\not\in\tau(A,Q)$. Then $q$ could be
  removed from the nucleus, contradicting its minimality.
\end{proof}

\begin{lem}[\cite{nekrashevych:ssg}*{Lemma~2.11.12}]
  If $\Pi$ is nuclear, then there exists constants $\lambda<1$, $K$
  and $n$ such that $|\tau(x,g)|<\lambda|g|+K$ for all $g\in\langle
  \Pi\rangle$ and $x\in A^n$, where $|g|$ denotes the minimal length
  of $g\in\langle \Pi\rangle$ as a word over $Q$.
\end{lem}

We introduce also the following
\begin{defn}
  Let $\Pi$ be a contracting automaton. Its set of states $Q=N$
  therefore contains a specific state, the identity, written
  $\varepsilon$. The automaton $\Pi$ is \emph{smooth} if $\tau$ is
  surjective and the following subgraph of $\gf(\Pi^*)$ is strongly
  connected: its vertices are letters $a\in A$; there is an edge from
  $a$ to $b$ for all $q\in Q$ with $\sigma(q,a)=b$ and
  $\tau(q,a)=\varepsilon$.
\end{defn}

The following consequence of smoothness will not be used directly; its
proof is implicit within the proof of
Proposition~\ref{prop:connected}: if an automaton is smooth, then it
is recurrent and spherically transitive. More is true:
\begin{lem}
  If $\Pi$ is a smooth automaton, then for every pair of
  alphabet-words $a,b\in A^*$ of the same length, and any word $v\in
  Q^*$, there exists a word $w\in Q^*$ with $\sigma(a,w)=b$ and
  $\tau(a,w)=\varepsilon^{n_0}v_1\varepsilon^{n_1}\dots
  v_k\varepsilon^{n_k}$ for some $n_0,\dots,n_k\ge0$.
\end{lem}

\subsection{Profinite groups}
Let a group $\Gamma$ act on a tree $A^*$. There is then a family of
finite quotients $\Gamma_n$ defined by restricting the action of
$\Gamma$ to $A^n$. These form a projective system
$\Gamma_{n+1}\to\Gamma_n$, with inverse limit
$\overline\Gamma=\varprojlim\Gamma_n$.

\begin{lem}\label{lem:closure}
  $\overline\Gamma$ is the closure of $\Gamma$ in the topological
  group $\aut(A^*)$, with its standard (compact-open) topology.
\end{lem}
\begin{proof}
  In the standard topology of $W=\aut(A^*)$, since $A^*$ is discrete,
  a basis of neighbourhoods of the identity is given by subgroups
  fixing larger and larger finite sets of vertices.  One may take as
  such subgroups the pointwise fixators $\stab_W(n)$ of $A^n$; and
  then $\Gamma/(\Gamma\cap \stab_W(n))\cong\Gamma_n$.
  
  $\overline\Gamma$ is then closed in $W$ because it is an inverse
  limit of the closed subgroups $\Gamma_n\subset W/\stab_W(n)$.
  
  To show that $\Gamma$ is dense in $\overline\Gamma$, pick
  $g=\varprojlim\gamma_n\in\overline\Gamma$, with
  $\gamma_n\in\Gamma_n$. Choose lifts $g_n\in\Gamma$ of $\gamma_n$.
  Since the action of $g_n$ agrees with that of $g$ on $A^n$, we have
  $g_n\to g$.
\end{proof}

%%%%%%%%%%%%%%%%%%%%%%%%%%%%%%%%%%%%%%%%%%%%%%%%%%%%%%%%%%%%%%%%
\section{Groupoids}
In what follows we will use interchangeably the words orbispace and
groupoid; more on that philosophy can be found
in~\cite{moerdijk:orbifolds}. Recall that a groupoid is a graph
$(X,G,s,t)$ with a multiplication $\setsuch{(g,g')\in G\times
  G}{s(g')=t(g)}\to G$, an inverse $(\;)^{-1}:G\to G$ and an identity
$1:X\to G$, such that
\[s(gh)=s(g),\quad t(gh)=t(h),\quad gg^{-1}=1_{s(g)},\quad
g^{-1}g=1_{t(g)},\quad 1_{s(g)}g=g=g1_{t(g)}.
\]
More concisely, a groupoid is a small category in which all arrows are
invertible. We will usually denote the groupoid simply by $G$.

The group $G_x=\setsuch{g\in G}{s(g)=t(g)=x}$ of self arrows of an
object $x\in X$ is called the \emph{isotropy} of $x$, and the subset
$G^x=\setsuch{g\in G}{s(g)=x}$ is called the \emph{fiber} of $x$.

By groupoid, we shall always mean topological groupoid, i.e.\ we shall
assume the sets $G,X$ will have natural topologies that make all
structural maps --- source, target, multiplication, inverse, unit ---
continuous.

As a trivial example of groupoid, we may take $G=X$ a topological
space, with $s,t$ the identity. An example that will be relevant to
the sequel is the following: consider
\begin{equation}\label{eq:I}
  X=\{0,1\}^\N,\quad I=\{(x,x)\}_{x\in X}\cup\{(w0\overline1,w1\overline0)\}_{w\in
  X^*}\cup\{(w1\overline0,w0\overline1)\}_{w\in X^*},
\end{equation}
where $\overline0$ and $\overline1$ denote the infinite words
$00\dots$ and $11\dots$ respectively. The set $X$ is given the
Tychonoff (product) topology, and $I$ is given the topology inherited
from $X\times X$. Note that this makes $I$ compact. Define $s(g)$ and
$t(g)$ as the projections on the first and second coordinate respectively;
set $(x,y)^{-1}=(y,x)$, and $(x,y)(y,z)=(x,z)$. We shall see
shortly~(Lemma~\ref{lem:interval}) that $I$ should be understood as
the interval $[0,1]$.

A morphism of groupoids $G\to G'$ is a pair of continuous maps
$\phi:G\to G'$, $\psi:X\to X'$, such that $\psi(s(g))=s(\phi(g))$,
$\psi(t(g))=t(\phi(g))$, $\phi(g)^{-1}=\phi(g^{-1})$,
and $\phi(gh)=\phi(g)\phi(h)$.

%Our groupoids will satisfy the further condition that $s$ (and
%therefore $t$) is locally a closed inclusion, namely that every $g\in
%G$ has a closed neighbourhood $U_g\ni g$ such that $s_{|U_g}$ is a
%homeomorphism onto a closed subset of $X$.

%\begin{lem}
%  In our groupoids, all fibers are discrete.
%\end{lem}
%\begin{proof}
%  For $g\in G^x$, consider $V$ an open subset of $U_g$ containing
%  $g$, where $U_g$ is chosen as above. Since $s_V$ is injective,
%  $G^x\cap V=\{g\}$, so all points in $G^x$ are open, and therefore
%  $G^x$ is discrete.
%\end{proof}

One can think of our groupoid (even though this is not equivalent) as
a topological space $Z=X/\{s(g)=t(g)\}$, together with for every $z\in
Z$ a small enough neighbourhood $U_z\ni z$ and a homeomorphism
$U_z\cong \tilde U_z/G_z$, where $G_z$ is a discrete group acting on a
topological space $\tilde U_z$ and fixing a point $\tilde z\in\tilde
U_z$ identified with $z$.

Let us give a vague outline of this construction in the special case
when $s$ is locally a homeomorphism. We denote by $\pi:X\to Z$ the
projection. For a point $z=\pi(x)\in Z$, we set $G_z=G_x$, and choose
$\tilde U_z$ as some neighbourhood of $\tilde z=x\in X$. We will
define an action of $G_z$ on this neighbourhood. Let $g\in G_z$ be an
arrow. Since $s:G\to X$ is a local homeomorphism with $s(g)=x$, it
admits a unique germ of section $\sigma_s:X\dashrightarrow G$, defined
on a neighbourhood of $x$, and satisfying $\sigma_s(x)=g$. We define
the action of $g$ to be the composite
$\phi_g=t\sigma_s:X\dashrightarrow X$. It is a germ of homeomorphism
defined on a neighbourhood of $x$ and fixing $x$. Clearly, if
$\phi_g(x')=x''$ then $s(g)=x'$ and $t(g)=x''$ and hence
$\pi(x')=\pi(x'')$, so $\pi$ is locally a homeomorphism $\tilde
U_z/G_z\to Z$ that\footnote{We are being sloppy here: first of all, it
  might not be possible to lift the action $G_z\to\{\text{germs of
    homeomorphisms}\}$ to an action defined on a neighbourhood of $x$.
  Secondly, $Z$ might fail to be Hausdorff, in which case we couldn't
  conclude that $\tilde U_z/G_z\to Z$ is injective.

The reader who is more concerned about mathematical rigor than
geometric interpretation may disregard the above construction at no
cost.} sends $\tilde z=x$ to $z=\pi(x)$.

\subsection{Geometric Realizations} The following classical
construction builds a topological space $|G|$ from a groupoid $G$,
such that all considerations on $G$ have equivalents on $|G|$. First,
let
\[G_n=\setsuch{(x_0,g_1,x_1,\dots,g_n,x_n)}{x_i\in X,g_i\in
  G,t(g_i)=x_i=s(g_{i+1})\text{ for all }i}\]
be the space of composable sequences of $n$ arrows. There are face maps
$d_i:G_n\to G_{n-1}$ and degeneracies $s_i:G_n\to G_{n+1}$, for
$i\in\{0,\dots,n\}$, given by
\begin{align*}
  d_0(x_0,g_1,x_1,\dots,g_n,x_n)&=(x_1,\dots,g_n,x_n),\\
  d_i(x_0,g_1,x_1,\dots,g_n,x_n)&=(x_0,\dots,x_{i-1},g_ig_{i+1},x_{i+1},\dots,g_n,x_n),\\
  d_n(x_0,g_1,x_1,\dots,g_n,x_n)&=(x_0,\dots,g_{n-1},x_{n-1}),\\
  s_i(x_0,g_1,x_1,\dots,g_n,x_n)&=(x_0,g_1,\dots,x_i,1_{x_i},x_i,\dots,g_n,x_n),\\
\end{align*}
which turn the family $(G_n)_{n\ge0}$ into a simplicial space. One
then lets $|G|$ be the geometric realization of that simplicial space;
namely, let $\Delta^n$ denote the standard $n$-simplex
$\Delta^n=\setsuch{(t_0,\dots,t_n)\in\R^{n+1}}{t_i\ge0,\sum t_i=1}$,
with its usual cofaces and codegeneracies
\begin{align*}
  \delta_i(t_0,\dots,t_{n-1})&=(t_0,\dots,t_{i-1},0,t_i,\dots,t_{n-1}),\\
  \sigma_i(t_0,\dots,t_{n+1})&=(t_0,\dots,t_{i-1}+t_i,\dots,t_{n+1});
\end{align*}
then we take
\[|G|=\coprod_{n\ge0}G_n\times
\Delta^n\bigg/\begin{array}{c}(d_i(x),t)\sim(x,\delta_i(t))\\
  (s_i(x),t)\sim(x,\sigma_i(t))\end{array}
\]
\begin{lem}
  If $G$ is compact metrizable, then so is the $k$-skeleton $|G|_k$ of
  $G$:
  \[|G|_k=\coprod_{n=0}^k G_n\times\Delta^n/{\sim}.\]
\end{lem}
\begin{proof}
  Let us assume by induction that $|G|_{k-1}$ is compact
  metrizable. The space $|G|_k$ is obtained from $|G|_{k-1}$ by
  gluing $G_k\times\Delta^k$ via an attaching map
  $f:G_k\times\partial\Delta^k\to|G|_{k-1}$. Therefore, $|G|_k$ is the
  pushout of the following diagram of compact spaces:
  \[\xymatrix{ & {G_k\times\Delta^k}\\
    {G_k\times\partial\Delta^k}\ar[ur]^{\iota}\ar[dr]_{f}\\
    & {|G|_{k-1}.}}
  \]
  Since $\iota$ is a closed inclusion, the pushout $|G|_k$ is the
  quotient of $|G|_{k-1}\sqcup G_k\times\Delta^k$ by a closed
  equivalence relation, and is thus compact and metrizable.
%%  The space $|G|_k$ is the quotient of the compact metrizable space
%%  $\coprod_{n=0}^k G_n\times\Delta^n$ by a closed relation:
%%  \[\coprod_{n=1}^k\coprod_{i=0}^nG_n\times\Delta^{n-1}\sqcup\coprod_{n=0}^{k-1}\coprod_{i=0}^nG_n\times\Delta^{n+1}\rightrightarrows\coprod_{n=0}^k
%%  G_n\times\Delta^n,\]
%%  where the two maps are
%%  $\alpha=\coprod_{n=1}^k\coprod_{i=0}^nd_i\times1_{\Delta^{n-1}}\sqcup\coprod_{n=0}^{k-1}\coprod_{i=0}^ns_i\times1_{\Delta^{n+1}}$ and 
%%  $\beta=\coprod_{n=1}^k\coprod_{i=0}^n1_{G_n}\times\delta_i\sqcup\coprod_{n=0}^{k-1}\coprod_{i=0}^n1_{G_n}\times\sigma_i$
%%  respectively.
\end{proof}

We say a groupoid $G$ is \emph{connected} if it is impossible to
disconnect $X=X'\sqcup X''$ and $G=G'\sqcup G''$ in disjoint open subsets
such that $(G',X')$ and $(G'',X'')$ are groupoids by restriction.
\begin{lem}
  $G$ is connected if and only if $|G|$ is connected.
\end{lem}
\begin{proof}
  If $|G|$ is not connected, write it as $|G|=A\sqcup B$, and set
  $X'=X\cap A$ and $X''=X\cap B$, where $X$ is viewed as a subspace of
  $|G|$. Set then $G'=\bigcup_{x\in X'}G^x$ and $G''=\bigcup_{x\in
  X''}G^x$. The target map sends $G'$ to $X'$ and $G''$ to $X''$,
  because otherwise there would be an edge crossing from $X'$ to
  $X''$, which is impossible by assumption.
  
  If $G$ is not connected, then $G=G'\sqcup G''$ and therefore
  $|G|=|G'|\sqcup|G''|$.
\end{proof}

\begin{defn}
  Let $(\phi,\psi)$ be a morphism of groupoids, $\phi:G\to G'$ and
  $\psi:X\to X'$. It is a \emph{covering} if $\psi$ is a covering, and
  \[\xymatrix{{G}\ar[r]^{\phi}\ar[d]_{s} & {G'}\ar[d]^{s'}\\
    {X}\ar[r]_{\psi} & {X'}}
  \]
  is a pull-back diagram.
\end{defn}
This implies that $\phi$ and $\psi$ have the same degree, which we
call the degree of $(\phi,\psi)$.  We shall actually abuse notation
and denote both $\phi$ and $\psi$ by the letter $\phi$.

\begin{lem}\label{lem:covering}
  Let
  \[\xymatrix{{\tilde R}\ar[d]_{p_R}\ar@<0.4ex>[r]^{\tilde\alpha}\ar@<-0.4ex>[r]_{\tilde\beta} &
    {\tilde X}\ar[d]^{p_X}\\
    {R}\ar@<0.4ex>[r]^{\alpha}\ar@<-0.4ex>[r]_{\beta} & {X}}
  \]
  be a diagram of compact spaces, where $\alpha p_R=p_X\tilde\alpha$,
  $\beta p_R=p_X\tilde\beta$, and $p_R,p_X$ are degree-$d$ covering
  maps. Suppose that $Z=X/R$ is Hausdorff, and that the map of
  coequalizers $p_Z:\tilde Z=\tilde X/\tilde R\to Z$ has fibers of
  cardinality $d$. Then $p_Z$ is a covering map.
\end{lem}
\begin{proof}
  Pick $z\in Z$ and let $C$ be its preimage in $X$. The restriction of
  $\tilde X$ to $C$ is trivialized by its map to $p_Z^{-1}(z)$. Extend
  that trivialization to an open neighbourhood $\U\supset C$.  It
  induces trivializations of $\tilde R$ on $\alpha^{-1}(\U)$ and
  $\beta^{-1}(\U)$. Moreover, these trivializations agree on
  $\alpha^{-1}(C)=\beta^{-1}(C)$. The structure group being finite
  hence discrete, they must agree on some neighbourhood
  $\V\supset\alpha^{-1}(C)$ contained in
  $\alpha^{-1}(\U)\cap\beta^{-1}(\U)$. Pick a neighbourhood $\U''$ of
  $C$ satisfying $\alpha^{-1}(\U'')\subset\V$ and
  $\beta^{-1}(\U'')\subset\V$.
  
  Next, saturate $\U''$ under the equivalence relation generated by
  $R$, by setting $\U'=X\setminus p_1p_2^{-1}(X\setminus\U'')$, where
  $p_1,p_2$ are respectively the first- and second-coordinate
  projections $X\times_Z X\rightrightarrows X$; they are closed maps
  on the compact $X\times_Z X$, because $Z$ was assumed Hausdorff.
 
  Set $\V'=\alpha^{-1}(\U')=\beta^{-1}(\U')$.
  The trivializations induce isomorphisms
  \begin{equation}\label{eq:lem:covering}
    \left(\raisebox{5ex}{\xymatrix{{\tilde R_{|\V'}}\ar[d]_{p_R}\ar@<0.4ex>[r]^{\tilde\alpha}\ar@<-0.4ex>[r]_{\tilde\beta} &
        {\tilde X_{|\U'}}\ar[d]^{p_X}\\
        {\V'}\ar@<0.4ex>[r]^{\alpha}\ar@<-0.4ex>[r]_{\beta} & {\U'}}}\right)\cong
    \left(\raisebox{5ex}{\xymatrix{{\V'\times\{1,\dots,d\}}\ar[d]\ar@<0.4ex>[r]^{\alpha\times1}\ar@<-0.4ex>[r]_{\beta\times1} &
        {\U'\times\{1,\dots,d\}}\ar[d]\\
        {\V'}\ar@<0.4ex>[r]^{\alpha}\ar@<-0.4ex>[r]_{\beta} & {\U'}}}\right).
  \end{equation}
  Therefore the coequalizers of~\eqref{eq:lem:covering} form a trivial
  covering. Now remember that the quotient $\U'/\V'$
  is a neighbourhood of $z\in Z$. We have just exhibited a
  trivialization of $\tilde Z$ on $\U'/\V'$. It
  follows that the map $p_Z:\tilde Z\to Z$ is a covering.
\end{proof}

\begin{prop}\label{prop:extend}
  A homomorphism $\phi:G\to G'$ induces functorially a continuous map
  $|\phi|:|G|\to|G'|$. The map $\phi$ is a covering if and only if
  $|\phi|$ is a covering, and then both have the same degree.
\end{prop}
\begin{proof}
  Let $(\phi:G\to G',\phi:X\to X')$ be our groupoid homomorphism. Define
  $|\phi|:|G|\to|G'|$ by
  \[|\phi|(x_0,g_1,\dots,x_n;t_0,\dots,t_n)=(\phi(x_0),\phi(g_1),\dots,\phi(x_n);t_0,\dots,t_n).\]
  
  Assume now that $\phi$ is a covering of degree some cardinal
  $\aleph$.  First note that composable sequences of arrows
  $(x_0',g_1',\dots,x_n')\in G_n'$ satisfy the ``unique lifting
  property''. Namely, given $x_0\in X$ with $\phi(x_0)=x_0'$ there is
  a unique composable sequence of arrows $(x_0,g_1,\dots,x_n)\in G_n$
  with $\phi(x_i)=x_i'$ and $\phi(g_i)=g_i'$.  Indeed
  $g_1'\in(G')^{\phi(x_0)}$ has a unique lift $g_1\in G^{x_0}$; we set
  $x_1=t(g_1)$, and lift $g_2'$, etc., until $x_n'$. Every sequence
  $(x_0',g_1',\dots,x_n')\in G_n'$ therefore has as many preimages in
  $G_n$ as $x_0'$ has preimages in $X$.
  
%%  We claim that the map $G_n\to G'_n$ is a covering of degree
%%  $\aleph$. To show that, it suffices to exhibit neighbourhoods
%%  $U_{(x_0,g_1,\dots,x_n)}$ of $(x_0,g_1,\dots,x_n)$ that map
%%  homeomorphically onto neighbourhoods $U_{(x'_0,g_1,\dots,x'_n)}$.
%%  Now by assumption each $\bullet\in\{x_0,g_1,\dots,x_n\}$ has a
%%  neighbourhood $U_\bullet$ that maps homeomorphically onto a
%%  neighbourhood $U_{\bullet'}$ of $\bullet'$, and we may set
%%  $U_{(x'_0,g_1,\dots,x'_n)}=(U_{x_0}\times U_{g_1}\dots\times
%%  U_{x_n})\cap G_n$.
  We have shown that $(G_n\to X)\to(G_n'\to X')$ is a pullback
  diagram, and hence that $G_n\to G_n'$ is a
  covering.\footnote{Actually, we have only shown it at the level of
    sets, but an easy diagram chase shows that it is also a pullback
    of topological spaces.}  Therefore, by crossing with $\Delta^n$,
  the map
  \[\coprod_{n\ge0}G_n\times\Delta^n\to\coprod_{n\ge0}G_n'\times\Delta^n\]
  is a covering of degree $\aleph$.

  Now recall that $|G|$ and $|G'|$ are obtained by quotienting, i.e.\
  taking a coequalizer
  \[\xymatrix{{\coprod_{n\ge1}\coprod_{i=0}^nG_n\times\Delta^{n-1}\sqcup\coprod_{n\ge0}\coprod_{i=0}^nG_n\times\Delta^{n+1}}\ar@<0.4ex>[r]^-{\alpha}\ar@<-0.4ex>[r]_-{\beta}
    & {\coprod_{n\ge0}G_n\times\Delta^n,}}
  \]
  where
  $\alpha=\coprod_{n=1}^k\coprod_{i=0}^nd_i\times1_{\Delta^{n-1}}\sqcup\coprod_{n=0}^{k-1}\coprod_{i=0}^ns_i\times1_{\Delta^{n+1}}$
  and
  $\beta=\coprod_{n=1}^k\coprod_{i=0}^n1_{G_n}\times\delta_i\sqcup\coprod_{n=0}^{k-1}\coprod_{i=0}^n1_{G_n}\times\sigma_i$.
  In the following diagram both rows are coequalizers:
  \[\xymatrix{{\displaystyle\coprod_{n=1}^k\coprod_{i=0}^nG_n\times\Delta^{n-1}\sqcup\coprod_{n=0}^{k-1}\coprod_{i=0}^nG_n\times\Delta^{n+1}}\ar@<0.4ex>[r]^(0.7){\alpha}\ar@<-0.4ex>[r]_(0.7){\beta}\ar[d]^{\dot\phi} &
    {\displaystyle\coprod_{n\ge0}G_n\times\Delta^n}\ar[r]\ar[d]^{\ddot\phi} &
    {\displaystyle|G|_k}\ar[d]^{|\phi|}\\
    {\displaystyle\coprod_{n=1}^k\coprod_{i=0}^nG'_n\times\Delta^{n-1}\sqcup\coprod_{n=0}^{k-1}\coprod_{i=0}^nG'_n\times\Delta^{n+1}}\ar@<0.4ex>[r]^(0.7){\alpha'}\ar@<-0.4ex>[r]_(0.7){\beta'} &
    {\displaystyle\coprod_{n\ge0}G'_n\times\Delta^n}\ar[r] &
    {\displaystyle|G'|_k}}
  \]
%%  The first two vertical arrows are coverings of degree $\aleph$, and
%%  the commuting squares $\alpha\ddot\phi=\dot\phi\alpha'$ and
%%  $\beta\ddot\phi=\dot\phi\beta'$ induce bijections on their
%%  fibers. It follows that the last vertical map $|\phi|$ is also a
%%  covering map, of same degree $\aleph$.
  The space $|G|_k$ is the union of the $G_n\times(\Delta^n)^\circ$.
  On each such piece, $|\phi|$ is the product of $\phi:G_n\to G_n'$
  and the identity on $(\Delta^n)^\circ$. In particular, $|\phi|$ is
  everywhere $d$-to-$1$. Applying Lemma~\ref{lem:covering}, we see
  that $|G|_k\to|G'|_k$ is a covering, and by taking direct limits, so
  is $|G|\to|G'|$.
  
  Conversely, if $|G|\to|G'|$ is a cover, then $X=|G|_0\to X'=|G'|_0$
  is a cover by restriction. The diagram $(G\to X)\to(G'\to X')$
  yields a map $G\to X\times_{X'}G'$, for which we need to find an
  inverse. The mapping cylinder of $s$ (respectively of $s'$)
  $M_s=G\times[0,1]\cup_s X$ sits as a subspace of $|G|$ (respectively
  $|G'|$). It is the image of $G\times [0,\frac12]\hookrightarrow
  G\times[0,1]\twoheadrightarrow|G|_1\subset|G|$. The map $M_s\to
  M_{s'}$ is therefore also a cover. We have a map
  $(X\times_{X'}G')\times[0,1]\to M_{s'}$ defined by
  $(x,g',t)\mapsto(g',t)$. Using the unique homotopy lifting
  extension, we obtain a lift $(X\times_{X'}G')\times[0,1]\to M_s$
  making the diagram
  \[\xymatrix{{(X\times_{X'}G')\times[0,1]}\ar@{.>}[r]\ar@/^2ex/[rr] & {M_s}\ar[r] & {M_{s'}}\\
    {(X\times_{X'}G')\times\{0\}}\ar[r]\ar@{^{(}->}[u] & {X}\ar[r]\ar@{^{(}->}[u] & {X'}\ar@{^{(}->}[u]}
  \]
  commute. Restriction to $(X\times_{X'}G')\to G\times\{1\}$ gives us
  our desired map.
\end{proof}

\begin{lem}\label{lem:interval}
  The groupoid $I$ defined in~\eqref{eq:I} admits an injective
  continuous map $\phi:[0,1]\to|I|$ with $\phi(0)=\overline0$ and
  $\phi(1)=\overline1$.
\end{lem}
\begin{proof}
  We first define a map $\chi:X\to[0,1]$, by
  \[\chi(w_1w_2\dots)=\sum_{i=1}^\infty2w_i3^{-i}.\]
  Then we define $\psi:|I|\to[0,1]$ by
  \[\psi(x_0,g_1,x_1,\dots,g_n,x_n;t_0,\dots,t_n)=\sum_{i=0}^nt_i\chi(x_i).\]
  Consider the image $|I|^+$ of $(\{(x,x)\}_{x\in
    X}\cup\{(w0\overline1,w1\overline0)\}_{w\in\{0,1\}^*})\times\Delta^1$
  in $|I|$; then $\psi$ defines a continuous bijection between $|I|^+$
  and $[0,1]$. Since $|I|^+$ is compact, $\psi$ admits a
  continuous inverse $\phi$.
\end{proof}

Actually, $\phi([0,1])$ is a deformation retract of $|I|$: the
$n$-skeleton $|I|_n$ is homeomorphic to the closure in $\R^{n+1}$ of
the countable union of spheres
\[\bigcup_{w\in\{0,1\}^*}\setsuch{x\in\R^{n+1}}{\rule{0pt}{2ex}\Big\|x-\big(\chi(w0\overline1)+\frac{3^{-|w|}}{6},0,\dots,0\big)\Big\|=\frac{3^{-|w|}}{6}}.
\]
%%\[\bigcup_{w\in\{0,1\}^*}\\mathbb S\left(w\overline
%%  0,\frac{3^{-|w|}}{6}\right),\text{ where }\mathbb
%%S(w,\rho)=\setsuch{x\in\R^{n+1}}{\rule{0pt}{2ex}\|x-(\chi(w)+\rho,0,\dots,0)\|=\rho}.
%%\]
\begin{center}
  \psfrag{0}{\large $0$}
  \psfrag{1}{\large $1$}
  \includegraphics{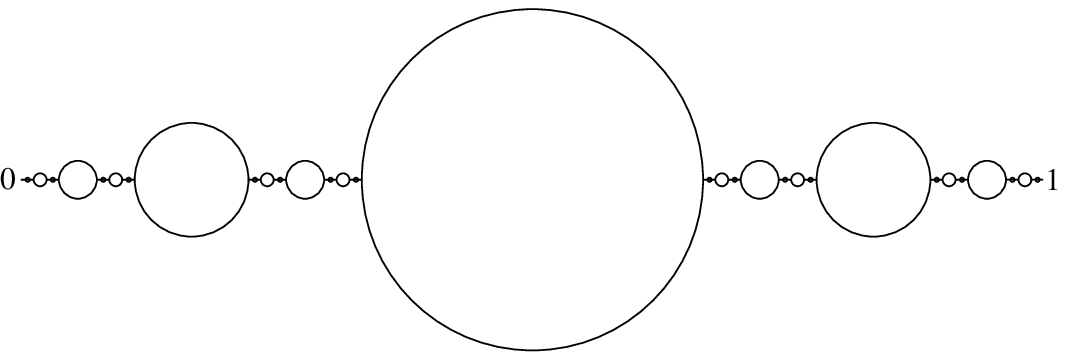}
\end{center}
A non-degenerate $n$-simplex $|\Sigma|$ of $|I|_n$ comes from a
sequence $\Sigma=(w0\overline1,w1\overline0,w0\overline1,\dots)$ or
$(w1\overline0,w0\overline1,w1\overline0,\dots)$ in $I_n$. There are
exactly two such simplices for every $(n-1)$-sphere in $I_{n-1}$,
which glue on $\Sigma$ as the two hemispheres attach to $|\Sigma|$ in
$|I|_n$. One can retract $|I|$ to $|I|^+_1$ by successively sliding
the $n$-spheres over the $(n+1)$-spheres, in shorter and shorter
amounts of time as $n\to\infty$ so that the retraction be performed in
finite total time.

We say that a groupoid $G$ is \emph{arcwise connected} if for every
$x,y\in X$ there exists a homomorphism $\gamma:I\to G$ with
$\gamma(\overline0)=x$ and $\gamma(\overline 1)=y$.
\begin{lem}\label{lem:arccon}
  $G$ is arcwise connected if and only if $|G|$ is arcwise connected.
\end{lem}
\begin{proof}
  Assume first that $G$ is arcwise connected.  Notice first that any
  point $(x_0,g_1,x_1,\dots,g_n,x_n;t_0,\dots,t_n)\in|G|$ is connected
  by an interval, within a simplex, to a point
  $(x_0,g_1,x_1,\dots,g_n,x_n;1,0,\dots,0)\sim(x_0,1)$ sitting in the
  $0$-skeleton $G_0\times\Delta^0=X$ of $G$. It is therefore enough to
  show that any two points $x,y\in X\subseteq|G|$ are connected by an
  interval.
  
  Now by assumption we have a homomorphism $\gamma:I\to G$ with
  $\gamma(\overline0)=x$ and $\gamma(\overline1)=y$, which induces a
  path $|\gamma|:|I|\to|G|$. Precomposing with the
  $\phi:[0,1]\to|I|$ from Lemma~\ref{lem:interval}, we obtain a path
  $\delta=|\gamma|\phi:[0,1]\to|G|$ with $\delta(0)=x$ and
  $\delta(1)=y$.
  
  Conversely, assume that $|G|$ is arcwise connected; let
  $x,y\in|G|_0$ be two points, and $\gamma:[0,1]\to|G|$ a path with
  $\gamma(0)=x$ and $\gamma(1)=y$. Set $J_n=\gamma^{-1}(|G|_n)$ for
  all $n\in\N$. Since the $J_n$ are closed in $[0,1]$ there is by
  compactness of $[0,1]$ a minimal $n\in\N$ with $J_n=[0,1]$. If
  $n>1$, we may perturb $\gamma$ on $J_n\setminus J_{n-1}$ to make it
  avoid
  $G_n\times\{(\frac1{n+1},\dots,\frac1{n+1})\}\setminus|G|_{n-1}$,
  and then homotope it within $\Delta^n$ so that its image is entirely
  contained in $|G|_{n-1}$. By repeating this argument, we may assume
  that $\gamma[0,1]$ lies in $|G|_1$ and is transverse to
  $G_1\times\{(\frac12,\frac12)\}\setminus|G|_0$.
  
  Set now $C=\gamma^{-1}(G_1\times\{(\frac12,\frac12)\})\setminus
  J_0$. It is discrete in $[0,1]\setminus J_0$, and therefore
  countable. Let $Y$ be the topological space obtained by cutting
  $[0,1]$ at every point in $C$; more precisely,
  $Y=[0,1]\times\{0\}\cup C\times\{1\}$ with the topology inherited
  from the lexicographical ordering. Let $\pi:Y\to[0,1]$ be projection
  on the first coordinate; then $\gamma\pi:Y\to|G|_1$ admits an
  obvious retraction to $\delta:Y\to|G|_0$. We let $K$ be the
  equivalence groupoid $Y/\pi$, namely
  $K=\setsuch{((t,\varepsilon),(t,\varepsilon'))}{(t,\varepsilon),(t,\varepsilon')\in
    Y}$; then the map $\delta$ extends to a morphism of groupoids
  $K\to G$ such that $\delta(0,0)=x$ and $\delta(1,0)=y$. It therefore
  suffices to exhibit a morphism $\iota$ of groupoids from $I$ to $K$
  with $\iota(\overline 0)=(0,0)$ and $\iota(\overline 1)=(1,0)$.
  
  Choose a surjective monotone map from the Cantor set $\{0,1\}^\N$ to
  $[0,1]$ such that all points in $C$ are hit by a point of the form
  $w0\overline1$ (and hence also of the form $w1\overline0$). This map
  lifts to a map $\{0,1\}^\N\to K$ and induces the desired $\iota$.
\end{proof}

\subsection{Topological quotients} Another way of associating a space
to a groupoid $(X,G)$ is to quotient $X$ by the equivalence relation
$x\sim y$ if there exists $g\in G$ with $s(g)=x$ and $t(g)=y$. The
\emph{topological quotient} of $G$ is $X/{\sim}$, denoted $G_\top$. It
is less well behaved than the geometric realization $|G|$, for example
because if $\phi:G\to H$ is a covering of orbispaces then $G_\top\to
H_\top$ is not necessarily a covering. It will typically be a
\emph{branched covering}, with branching locus the set of points $x\in
G_\top$ where $G_x\to H_{\phi(x)}$ is not a group
isomorphism. However, we still have the equivalences
\begin{prop}
  \begin{gather*}
    G\text{ connected }\Leftrightarrow|G|\text{ connected
    }\Leftrightarrow G_\top\text{ connected,}\\
    G\text{ compact }\Leftrightarrow|G|_1\text{ compact
    }\Leftrightarrow|G|_k\text{ compact }\forall k\in\N\Rightarrow
    G_\top\text{ compact.}
  \end{gather*}
\end{prop}
Note that $\pi_1(|G|)$ and $\pi_1(G_\top)$ are very different; we are
more interested in the former, and denote $\pi_1(G)=\pi_1(|G|)$.

\subsection{Morita Equivalence}
In category theory, two categories can be equivalent even if their
objects cannot be put in bijection.\footnote{For instance, the trivial
  category with one object and one isomorphism is equivalent to the
  category with two objects $a,b$ and four isomorphisms, one between
  any two objects.} The corresponding notion for topological groupoids
is usually called \emph{Morita equivalence}. Unlike abstract
categories, two Morita equivalent groupoids do not always admit
continuous functors between them:
\begin{defn}
  Let $(G,X)$ and $(G',X')$ be two topological groupoids. They are
  \emph{Morita equivalent} if there exists a topological space $P$,
  equipped with two maps $s_P:P\to X$ and $t_P:P\to X'$, such that
 \[\xymatrix{G''=G\sqcup P\sqcup P\sqcup
  G'\ar@<1ex>[d]^{s_G\sqcup s_P\sqcup t_P\sqcup s_{G'}}\ar@<-1ex>[d]_{t_G\sqcup t_P\sqcup
  s_P\sqcup t_{G'}}\\ X''=X\sqcup X'}
  \]
  can be endowed with the structure of a groupoid. We further assume
  that the bijective maps $P/(p\sim ph)\to X$ and $P/(p\sim gp)\to X'$
  are homeomorphisms.
  
  Assume furthermore that $G$ and $G'$ come equipped with
  self-coverings $f,f'$. Then $(G,f)$ and $(G',f')$ are
  \emph{Morita equivalent} if $(G'',X'')$ admits a self-covering
  $f''$ whose restriction to $G$ and $G'$ yields $f$ and $f'$
  respectively.
\end{defn}

For instance, if $X=X'=\{\cdot\}$, then $G$ and $G'$ are Morita
equivalent precisely when they are isomorphic groups. As another
example, the groupoid $I$ defined in~\eqref{eq:I} is Morita equivalent
to the groupoid $G'=X'=[0,1]$, by taking $P=\{0,1\}^\N$ with
$s_P(w)=w$ and $t_P(w)=\sum w_i2^{-i}$.

\begin{lem}
  If $G$ and $G'$ are two Morita equivalent groupoids, then
  \begin{enumerate}
  \item their topological quotients $G_\top$ and $(G')_\top$ are
    homeomorphic;
  \item there is a functorial bijection $\mathcal J:\{\text{covers of
    }G\}\leftrightarrow\{\text{covers of }G'\}$;
  \item if furthermore $(G,f)$ and $(G',f')$ are Morita equivalent,
    then $\mathcal J$ associates the cover corresponding to $f^n$ with
    the cover corresponding to $(f')^n$ for all $n\in\N$.
  \end{enumerate}
\end{lem}
\begin{proof}
%%  We view points of $G_\top$ as $G$-orbits in $X$. For $\orbit$ a
%%  $G$-orbit, we set $\phi(\orbit)=t_P(s_P^{-1}(\orbit))$. This
%%  is a bijection, because $t_P$ and $s_P$ are onto. The closed sets
%%  $F$ of $G_\top$ are $G$-invariant closed sets of $X$, and are mapped
%%  to $G'$-invariant sets of $X'$ by $F\mapsto t_P(s_P^{-1}(F))$. These
%%  are closed because $t_P$ is closed, so $\phi$ is closed. The same
%%  argument starting from $G'$ shows that $\phi^{-1}$ is closed.
  For the first point, it suffices to write
  \begin{align*}
    G_\top &= X/(s(g)\sim t(g))=(P/(p\sim ph))/(p\sim gp)\\
    (G')_\top &= X'/(s(h)\sim t(h))=(P/(p\sim gp))/(p\sim ph).
  \end{align*}
  
  For the second point, let $(H,Y)$ be a cover of $(G,X)$ with maps
  $\pi:H\to G$ and $\pi:Y\to X$. Define
  \begin{align*}
    Y'&=\setsuch{(y,p)\in Y\times P}{s_P(p)=\pi(y)}\big/\forall h\in
    H:\,(t(h),p)\sim(s(h),\pi(h)p),\\
    H'&=\setsuch{(h,p,q)\in H\times P^2}{s_P(p)=s(\pi h),\,s_P(q)=t(\pi
      h)}\big/\\
    &\phantom{\setsuch{(h,p,q)\in H\times P^2}{}}
    \forall k,\ell\in H:\,(h,p,q)\sim(kh\ell^{-1},\pi(k)p,\pi(\ell)q),
  \end{align*}
  which is a groupoid, with source $s(h,p,q)=(s(h),p)$, target
  $t(h,p,q)=(t(h),q)$, inverse $(h,p,q)^{-1}=(h^{-1},q,p)$ and
  multiplication $(h,p,q)\cdot(k,q,r)=(hk,p,r)$. This construction
  associates to $(H,Y)$ a groupoid $(H',Y')$, with a cover
  $\pi':(H',Y')\to(G',X')$ given by $\pi'(y,p)=t_P(p)\in X'$ and
  $\pi'(h,p,q)=p^{-1}\pi(h)q\in G'$.
  
  We now show that the map $\pi':Y'\to X'$ is a cover. Indeed, it is the
  coequalizer of the diagram
  \[\xymatrix{{H\times_XP}\ar@<0.4ex>[r]\ar@<-0.4ex>[r]\ar[d] &
    {Y\times_XP}\ar[r]\ar[d] & {Y'}\ar[d]\\
    {G\times_XP}\ar@<0.4ex>[r]\ar@<-0.4ex>[r] & {P}\ar[r]^{t_P} & {X'}}
  \]
  and the first two vertical arrows are pulled back from $Y\to X$,
  which was assumed to be a cover.
  
  To apply Lemma~\ref{lem:covering}, we need to show that $Y'\to X'$
  is $d$-to-$1$. We first show that the fiber $(\pi')^{-1}(x')$ is
  isomorphic to the fiber $\pi^{-1}(s_P(p))$ via $F:y\mapsto[(y,p)]$,
  for some $p\in P$ satisfying $t_P(p)=x'$. To see that $F$ is onto,
  notice that any $(\overline y,\overline p)\in(\pi')^{-1}(x')$ is
  equivalent to $(s(h),p)$ for some $h$ lifting $p\overline p^{-1}$
  and having target $y$. Also, the relation
  $(t(h),p)\sim(s(h),\pi(h)p)$ is transitive, and the only way to have
  $(y,p)\sim(y',p)$ is to let $h=1$, so $F$ is injective. It follows
  that $Y'\to X'$ is $d$-to-$1$ and by Lemma~\ref{lem:covering} it is
  a cover.

  We next show that this construction is a bijection. For this, apply
  it again to $(H',Y')$, yielding a groupoid $(H'',Y'')$. This new
  groupoid is canonically isomorphic to $(H,Y)$. We show this on the
  level of $Y''$, and skip the details for $H''$.
  
  By definition,
  \begin{multline}\label{eq:y''}
    Y''=\setsuch{(y,p,p')\in Y\times
      P^2}{s_P(p)=\pi(y),\,t_P(p)=t_P(p')}\big/\\
    \forall h\in H,g\in G':\,(t(h),p,p')\sim(s(h),\pi(h)pg,p'g).
  \end{multline}
  We defined two maps $Y\to Y''$, $Y''\to Y$ and show that they are
  inverses to each other. First, given $y\in Y$, there exists $p\in P$
  with $s_P(p)=\pi(y)$, and we define $Y\to Y''$ by $y\mapsto
  (y,p,p)$. If we had made another choice $p'\in P$, we would have
  obtained $(y,p',p')$ which is equivalent to $(y,p,p)$ by taking
  $h=1_y,g=p^{-1}p'$ in~\eqref{eq:y''}; so the map is well defined.
  
  Given $(y,p,p')\in Y''$, there exists a unique $k\in H$ with
  $s(k)=y$ and $\pi(k)=p(p')^{-1}$, because $\pi$ is a covering. We
  define $Y''\to Y$ by $(y,p,p')\mapsto t(k)$. If in~\eqref{eq:y''} we
  consider two equivalent elements $(t(h),p,p')$ and
  $(s(h),\pi(h)pg,p'g)$, and find $k\in H$ with $s(k)=t(h)$ and
  $\pi(k)=p(p')^{-1}$, then $s(hk)=s(h)$ and
  $\pi(hk)=\pi(h)pg(p'g)^{-1}$; then $t(hk)=t(k)$ so the map is well
  defined.
  
  Next, $Y\to Y''\to Y$ is the identity, because
  $y\mapsto(y,p,p)\mapsto y$ using $k=1_y$ in the paragraph above.
  
  Finally, $Y''\to Y\to Y''$ is the identity: we have $(y,p,p')\mapsto
  t(k)\mapsto(t(k),p',p')$ where $s(k)=y$ and $\pi(k)=p(p')^{-1}$. We
  take $h=k,g=1_{t(p')}$ to show that $(t(k),p',p')$ and $(y,p,p')$
  are equivalent.
  
  We now check the last assertion of the lemma. We again check it only
  on objects; the proof is similar for morphisms. Let $f''$ be the
  covering map of $P$ given by the Morita equivalence. The cover of
  $(G,X)$ associated to $f^n$ is again $(G,X)$, which yields
  \[Y'=\setsuch{(x,p)\in X\times P}{s_P(p)=f^n(x)}\big/\forall g\in
  G:\,(t(g),p)\sim(s(g),f^n(g)p),
  \]
  with covering map $Y'\to X'$ given by $(x,p)\mapsto t_P(p)$. We show
  that this covering is equivalent to $X'\to X'$, $x'\mapsto(f')^nx'$.
  
  Given $x'\in X'$, there exists $q\in P$ with $t_P(q)=x'$. Define
  $X'\to Y'$ by $x'\mapsto(s_P(q),(f'')^n(q))$.  Conversely, given
  $(x,p)\in Y'$, there exists a unique $q\in P$ with $(f'')^n(q)=p$
  and $s_P(q)=x$, because $f''$ is a cover. Define $Y'\to X'$ by
  $(x,p)\mapsto t_P(q)$. These two maps are mutual inverses.
\end{proof}

%%%%%%%%%%%%%%%%%%%%%%%%%%%%%%%%%%%%%%%%%%%%%%%%%%%%%%%%%%%%%%%%
\section{Orbispaces from automata}\label{sec:o<a}
Let $\Pi$ be a finite invertible automaton. The \emph{action space} of
$\Pi$ is the groupoid $O(\Pi)$ with objects $A^\N$, with
$\N=\{1,2,\dots\}$, with morphisms
\[O(\Pi)=\setsuch{(\alpha,\phi)\in A^\N\times
  \langle\Pi\rangle^\N}{\phi(n+1)=\tau(\alpha(n),\phi(n))\text{ for
    all }n},
\]
with source map $s(\alpha,\phi)=\alpha$ and target map
$t(\alpha,\phi)(n)=\sigma(\alpha(n),\phi(n))$.

The \emph{limit space} of $\Pi$ is the groupoid $L(\Pi)$ with objects
$A^{-\N}$, morphisms
\begin{multline*}
  L(\Pi)=\left\{(\alpha,\phi)\in A^{-\N}\times
    \langle\Pi\rangle^{-\N}\right|\phi(n+1)=\tau(\alpha(n),\phi(n))\text{
    for all $n$,}\\
  \text{ and $\phi(-\N)$ is finite}\big\}
\end{multline*}
with source map $s(\alpha,\phi)=\alpha$ and target map
$t(\alpha,\phi)(n)=\sigma(\alpha(n),\phi(n))$. Note
that by our convention $0\notin\N$. We extend $\phi$ at $0$ by
$\phi(0)=\tau(\alpha(-1),\phi(-1))$.

The \emph{solenoid space} of $\Pi$ is the groupoid $S(\Pi)$ with
objects $A^\Z$, morphisms
\begin{multline*}
  S(\Pi)=\left\{(\alpha,\phi)\in A^\Z\times
    \langle\Pi\rangle^\Z\right|\phi(n+1)=\tau(\alpha(n),\phi(n))\text{ for
    all $n$,}\\
  \text{ and $\phi(\Z)$ is finite}\big\},
\end{multline*}
with source map $s(\alpha,\phi)=\alpha$ and target map
$t(\alpha,\phi)(n)=\sigma(\alpha(n),\phi(n))$.

Note that the finiteness condition ``$\phi(\N)$ is finite'' is
automatically satisfied in $O(\Pi)$.

The sets $A$ and $\langle\Pi\rangle$ are given the discrete topology.
For all of these groupoids, the object set is given the Tychonoff
(product) topology, in which a basis for the topology is given by
cylinders
\[\orbit_{n_1,a_1,\dots,n_k,a_k}=\setsuch{\alpha:(\N\text{ or
  }-\N\text{ or }\Z)\to A}{\alpha(n_1)=a_1,\dots,\alpha(n_k)=a_k}.
\]
Similarly, the morphism set is given the restriction of the product
topology.

The action space is the usual groupoid considered in association with a
group action: there is a morphism from $\alpha\in A^\N$ to $g(\alpha)$
for all $g\in\langle\Pi\rangle$. The action space has few useful
topological properties, because in most cases the orbits of
$\langle\Pi\rangle$ are dense in $A^\N$.

The limit space is endowed with a unilateral shift $f:L(\Pi)\to
L(\Pi)$, given by $f(\alpha,\phi)(n)=(\alpha(n-1),\phi(n-1))$. This
map defines a $\#A$-to-$1$ self-covering map of $L(\Pi)$.

The solenoid is endowed with a bilateral shift $f:S(\Pi)\to S(\Pi)$,
given by $f(\alpha,\phi)(n)=(\alpha(n-1),\phi(n-1))$. This map defines
a homeomorphism of $S(\Pi)$.

\begin{lem}[\cite{nekrashevych:ssg}*{Proposition~3.2.8}]
  If $\Pi$ is nuclear with nucleus $N$, then the groupoids
  $L(\Pi)$ and $S(\Pi)$ have finite fibers of cardinality at most
  $\#N$.
  
  Therefore the topological dimension of $L(\Pi)_\top$ and of
  $S(\Pi)_\top$ is at most $\#N-1$.
\end{lem}

%% Let $(G,X)$ be a locally compact groupoid. It is
%% \emph{amenable}~\cite{anatharaman-r:amenable} if there exists a
%% family $(\mu_n^x)_{n\in\N,x\in X}$ of probability measures on $G^x$,
%% with $x\mapsto\mu_n^x(f)$ continuous for all continuous, compactly
%% supported function $f\in C_c(G)$, with
%% \[\lim_{n\to\infty}\|g\mu_n^{s(g)}0\mu_n^{r(g)}\|=0\text{ uniformly on
%%   compact subsets of }G.
%% \]
%% \begin{prop}
%%   If $\Pi$ is nuclear, then $O(\Pi)$ is amenable.
%% \end{prop}
%% \begin{proof}
%%   By~\cite{bartholdi-g-n:fractal}*{Proposition~8.9} every schreier
%%   graph, i.e.\ leaf of $O(\Pi)$, has polynomial growth. Therefore
%%   $\mu_n^x$ may be defined ...
%% \end{proof}

\begin{lem}\label{lem:contract}
  If $\Pi$ is nuclear, then
  \begin{gather*}
    L(\Pi)=\setsuch{(\alpha,\phi)\in A^{-\N}\times
      Q^{-\N}}{\phi(n+1)=\tau(\alpha(n),\phi(n))\text{ for all $n$}},\\
    S(\Pi)=\setsuch{(\alpha,\phi)\in A^\Z\times
      Q^\Z}{\phi(n+1)=\tau(\alpha(n),\phi(n))\text{ for all $n$}}.
  \end{gather*}
\end{lem}
In other words, the morphisms of $L(\Pi)$ are given by
negative-infinite paths in $\gf(\Pi)$, and the morphisms of $S(\Pi)$
are given by bi-infinite paths in $\gf(\Pi)$. In the other direction,
the morphisms of $O(\Pi)$ are generated by positive-infinite paths in
$\gf(\Pi)$, since these correspond to generators of
$\langle\Pi\rangle$.

\begin{cor}
  If $\Pi$ is recurrent and nuclear, then $L(\Pi)$ is connected,
  metrizable, with topological dimension at most $\#Q-1$.
\end{cor}

\begin{thm}[\cite{nekrashevych:ssg}*{Proposition~5.7.8}]
  \[S(\Pi)=\varprojlim L(\Pi)\stackrel{T}{\longleftarrow}L(\Pi)\stackrel{T}{\longleftarrow}\cdots.\]
\end{thm}

%% Recall the following definition from~\cite{ruelle:thermo}, also
%% described in~\cite{putnam:smale}:
%% \begin{defn}
%%   A \emph{Smale space} is a compact metric space $X$ with a
%%   homeomorphism $\phi$ (one usually assumes furthermore that $\phi$ is
%%   topologically mixing), and constants $\epsilon_0>0$ and
%%   $\lambda_0\in(0,1)$, and a continuous map $[,]$:
%%   \[(x,y)\in\setsuch{(x,y)\in X\times
%%     X}{d(x,y)\le2\epsilon_0}\to[x,y]\in X\]
%%   satisfying the axioms
%%   \begin{gather*}
%%     [x,x]=x,\quad [[x,y],z]=[x,z]=[x,[y,z]]\text{ whenever defined},\\
%%     [\phi x,\phi y]=\phi[x,y]\text{ whenever defined},\\
%%     d(\phi y_1,\phi y_2)\le\lambda_0d(y_1,y_2)\text{ whenever
%%     }[x,y_i]=y_i\text{ and }d(x,y_i)<\epsilon_0,\\
%%     d(\phi^{-1} y_1,\phi^{-1} y_2)\le\lambda_0d(y_1,y_2)\text{ whenever
%%     }[y_i,x]=y_i\text{ and }d(x,y_i)<\epsilon_0.
%%   \end{gather*}
%% \end{defn}
%% \fobx{MORE EXPLANATIONS. RECALL STABLE, UNSTABLE COORDS}
%% 
%% \begin{thm}  
%%   If $\Pi$ is nuclear, then $S(\Pi)$ is Smale space.
%% \end{thm}

Consider the Schreier graphs $\gf_n$ of $\Gamma=\langle\Pi\rangle$
on $A^n$. There is a graph covering map $\gf_{n+1}\to\gf_n$, given by
the map
\[a_1\dots a_{n_+1}\mapsto a_1\dots a_n\]
on their vertices. There is also a map between the geometric
realizations\footnote{Here we mean the usual geometric realization of
  a graph, not to be confused with the geometric realization $|G|$
  of a groupoid $G$.} of $\gf_{n+1}$ and $\gf_n$, given as follows:
first, map the vertices by
\[a_1a_2\dots a_{n+1}\mapsto a_2\dots a_{n+1}.\]
Then, map the edge labeled $q:a_1\dots a_{n+1}\to b_1\dots b_{n+1}$
to the edge labeled $\tau(a_1,q):a_2\dots a_{n+1}\to b_2\dots
b_{n+1}$.

Assume now that $\Pi$ is smooth. Then an \emph{expansion rule} for
$\Pi$ is the following: first, for every $q\in Q$, elements $e_q\in A$
and $v_q\in Q$ with $\tau(e_q,v_q)=q$.  Second, for every $a,a'\in A$
a word $w_{a,a'}\in Q^*$ such that $\sigma(a,w_{a,a'})=a'$ and
$\tau(a,w_{a,a'})=\varepsilon^{|w_{a,a'}|}$, a power of the identity
state.

We note that the existence of an expansion rule is equivalent to the
smoothness of the automaton.

\begin{prop}\label{prop:schreierlim}
  If $\Pi$ is nuclear, then $L(\Pi)$, as a topological
  graph\footnote{i.e.\ as a topological space $X$ with edges
    $E\rightrightarrows X$}, is the inverse limit of the geometric
  realizations of the Schreier graphs $\gf_m$ of $\langle\Pi\rangle$
  on $A^m$.
\end{prop}
\begin{proof}
  Clearly $A^{-\N}=\varprojlim A^m$, where the map $A^{m+1}\to A^m$ is
  given by deletion of the first letter. Note that the edges of
  $\gf_m$ are in one-to-one correspondence with the set
  \begin{equation}\label{eq:schreierlim}
    \setsuch{(\alpha,\phi)\in A^m\times
      Q^m}{\phi(n+1)=\tau(\alpha(n),\phi(n))\text{ for all
        $n\in\{1,\dots,m-1\}$ }}
  \end{equation}
  where to any edge $q$ of $\gf_m$ from $a_1\dots a_m$ to $b_1\dots
  b_m$ one puts in correspondence the groupoid element $(\alpha,\phi)$
  given by $\alpha(n)=a_n$ and $\phi(n)=\tau(a_1,\dots,a_{n-1},q)$. It
  is then clear by Lemma~\ref{lem:contract} that the inverse limit
  of~\eqref{eq:schreierlim} is $L(\Pi)$.
\end{proof}

\begin{prop}\label{prop:connected}
  If $\Pi$ is nuclear and spherically transitive, then $L(\Pi)$ is
  connected.  If furthermore $\Pi$ is smooth, then $L(\Pi)$ is arcwise
  connected.
\end{prop}
\begin{proof}
  To prove that $L(\Pi)$ is connected it is sufficient, thanks to
  Proposition~\ref{prop:schreierlim}, to show that the graphs $\gf_n$ are
  connected; but this is precisely the condition that $\Pi$ is
  spherically transitive.
  
  Assume now that $\Pi$ is smooth, and let
  $(\{e_q\},\{v_q\},\{w_{a,a'}\})$ be an expansion rule for $\Pi$. By
  Lemma~\ref{lem:arccon}, it is sufficient to construct, for any
  points $x,y\in A^{-\N}\subset L(\Pi)$, a path $\gamma:[0,1]\to
  |L(\Pi)|$ from $x$ to $y$. We will define partial maps
  $\gamma_n:[0,1]\dasharrow|L(\Pi)|$ converging to $\gamma$.
  
  Each $\gamma_n$ will be defined on a finite union of closed
  subintervals of $[0,1]$, in such a way that if $[a,b]$ and $[c,d]$
  are two consecutive intervals, then $\gamma_n(b)$ and $\gamma_n(c)$
  are in the $0$-skeleton of $|L(\Pi)|$ and have identical last $n$
  symbols.
  
  We start by $\gamma_0$ defined only at $0$ and $1$, with
  $\gamma_0(0)=x$ and $\gamma_0(1)=y$. Assume now that $\gamma_{n-1}$
  has been defined; $\gamma_n$ coincides with $\gamma_{n-1}$ on the
  domain of $\gamma_{n-1}$. Consider two consecutive intervals $[a,b]$
  and $[c,d]$ on which $\gamma_{n-1}$ is defined, and write
  $\gamma_{n-1}(b)=\dots b_{-n}b_{-n+1}\dots b_{-1}$ and
  $\gamma_n(c)=\dots c_{-n}b_{-n+1}\dots b_{-1}$. Set
  $u=w_{b_{-n},c_{-n}}$ and $\ell=|u|$, and cut the interval $[b,c]$
  into $2\ell+1$ parts $E_0,F_1,E_1,\dots,F_\ell,E_\ell$. Define
  $\gamma_n$ on $F_i$ as the linear map from $F_i$ onto the geometric
  realization of the edge $(\alpha,\phi)$ defined by
  \[\alpha(m)=\begin{cases}b_m&\text{ if }m>-n,\\
    \sigma(b_{-n},u_1\dots u_{i-1})&\text{ if }m=-n,\\
    e_{\phi(m+1)}&\text{ if }m<-n,\end{cases}\qquad
  \phi(m)=\begin{cases}\varepsilon&\text{ if }m>-n,\\
    u_i&\text{ if }m=-n,\\
    v_{\phi(m+1)}&\text{ if }m<-n.\end{cases}\]
  
  There is clearly a partially defined map
  $\gamma'(t)=\lim_{n\to\infty}\gamma_n(t)$. On the intervals at which
  it is not defined, it can be extended by a constant path; we let
  $\gamma$ be this extension $[0,1]\to|L(\Pi)|$.
  
  Because of smoothness, for every interval $[b,c]$ between two
  intervals of definition of $\gamma_n$ there is a corresponding word
  $a\in A^n$ for which the following happens: Given a point
  $x\in[b,c]$, either $x$ is mapped to an element
  $\alpha\in|L(\Pi)|_0=A^{-\N}$ of the $0$-skeleton whose last $n$
  entries agree with $a$, or $x$ is mapped to some edge
  $(\alpha,\phi)\in|L(\Pi)|_1$, where the last $n$ entries of $\alpha$
  agree with $a$ and the last $n$ entries of $\phi$ are all the
  identity state $\varepsilon$.
  
  It remains to show that $\gamma$ is continuous; intuitively, this
  happens because, as $n\to\infty$ and the length of the intervals
  $F_i$ on which the $\gamma_n$ are defined tends to $0$, the images
  of the $\gamma_n$ are paths in $|L(\Pi)|$ which become closer and
  closer to identity morphisms, and degenerate to a single point in
  $|L(\Pi)|$.
  
  Let us be more precise. The map $\gamma:[0,1]\to|L(\Pi)|_1$ is
  clearly continuous on all $F_i$'s interiors, so there remains to
  consider two cases: first, continuity at a point
  $x\in[0,1]\setminus\bigcup F_i$.  Its image $\gamma(x)$ lies in
  $A^*=|L(\Pi)|_0$. Viewing $|L(\Pi)|_1$ as a quotient of
  $L(\Pi)_1\times\Delta^1$ by the relation $\sim$ that collapses all
  sets $\{1_w\}\times\Delta^1$ to a point, a neighbourhood of
  $\gamma(x)=1_{\gamma(x)}\times\Delta^1/{\sim}$ can be chosen of the
  form $V\times\Delta^1/{\sim}$, where $V$ is a neighbourhood of
  $1_{\gamma(x)}$ in $L(\Pi)_1$. We may assume furthermore $V$ to be
  the cylinder of all elements in $L(\Pi)_1\subset A^{-\N}\times
  Q^{-\N}$ with prescribed last $n$ terms: the set of $(\alpha,\phi)$
  in which the last $n$ letters of $\alpha$ agree with $w$, and the
  last $n$ states of $\phi$ are all $\varepsilon$. Then
  $\gamma^{-1}(V\times\Delta^1/{\sim})$ contains the interval
  $(b,c)\ni x$, where $b$ and $c$ are chosen such that $[a,b]$ and
  $[c,d]$ are two consecutive intervals on which $\gamma_n$ is defined
  and $b<x<c$.  We have shown that the preimage of a neighbourhood of
  $\gamma(x)$ is a neighbourhood of $x$, so $\gamma$ is continuous at
  $x$.
  
  Next, let us consider both left and right continuity at
  $x\in\partial F_i$. On the side of $F_i$ it is clear from the
  definition of the $\gamma_n$, and on the other side the above
  argument applies.
\end{proof}

Note that we have actually shown slightly more in this proof; namely,
since all the morphisms $(\alpha,\phi)$ considered satisfy
$\phi_0=\varepsilon$, we have shown, anticipating
Definition~\ref{defn:stdtile}:
\begin{prop}\label{prop:tileconn}
  If $\Pi$ is nuclear and smooth\footnote{Recall that smoothness
    implies spherical transitivity} then the standard tile
  $\tile(\Pi)$ is arcwise connected.
\end{prop}

\section{Tilings from automata}\label{sec:t<a}
This section extends on Section~\ref{sec:o<a}, where a limit orbispace
$L(\Pi)$ and a solenoid $S(\Pi)$ were constructed from an automaton
$\Pi$. First, we define a few more geometric objects associated with an
automaton $\Pi$.
\begin{defn}\label{defn:stdtile}
  The \emph{standard tile} $\tile(\Pi)$ is the groupoid with same
  objects $A^{-\N}$ as the limit space, but whose morphisms are
  inverse sequences starting at the trivial state:
  \[\tile(\Pi)=\setsuch{(\alpha,\phi)\in L(\Pi)}{\phi(0)=\varepsilon}.\]
\end{defn}

The group $\Gamma$ embeds in the group of homeomorphisms of $A^\N$.  The
\emph{germ space} $\germ(\Pi)$ has as objects $A^\N$, and has as set
of morphisms from $w$ to $w'$ all germs of homeomorphisms mapping $w$
to $w'$ and coming from $\Gamma$. It is a quotient of the action groupoid:
for instance, in the action groupoid, the automorphisms at $w$ is the
isotropy subgroup of $w$, while in the germ groupoid it is the
quotient of the isotropy of $w$ by the \emph{stable isotropy} of $w$,
i.e.\ those $g\in \Gamma$ that fix an open neighbourhood of $w$.

The geometric notions corresponding to the tiles are slightly more
tricky to set up: the translates of the standard tile almost cover the
solenoid, but for the discarded edges. We want to add ``half of each
edge'' to the standard tile to obtain a covering of the solenoid with
empty-interior intersections.

Call an edge $(\alpha,\phi)$ of $L(\Pi)$ \emph{critical} if
$\phi_0\neq\varepsilon$, and define the critical locus
$C\subset|L(\Pi)|$ as follows. First let $C_1\subset|L(\Pi)|_1$ be the
set of middle points of critical edges.  A point
$x=(g,t)\in|L(\Pi)|_n$ with $g\in L(\Pi)_n$ and $t\in\Delta^n$ is in
$C$ if (one of) the closest point(s) to $t$ in the $1$-skeleton of
$\Delta^n$, call it $t'$, satisfies $(g,t')\in C_1$.  The critical
locus is a hypersurface in $|L(\Pi)|$ that intersects all critical
edges in their middle points.

\begin{defn}
  The \emph{geometric standard tile} $\overline\tile$ is obtained
  from $|L(\Pi)|$ by cutting it along the critical locus $C$. We
  therefore have a surjective map $\overline\tile\to|L(\Pi)|$ that is
  a homeomorphism away from $C$ and that is generically 2-to-1 on $C$.
\end{defn}

\begin{lem}
  The natural inclusion $|\tile(\Pi)|\subset\overline\tile$ is a
  deformation retract.
\end{lem}
\begin{proof}
  Consider a simplex
  $K=(x_0,g_1,\dots,x_n)\times\Delta^n\subset|L(\Pi)|$. We will define
  a retraction on each connected component $K'$ of the preimage of $K$
  in $\overline\tile$. We view $K'$ as a subset of $K$. Let $A$ and
  $B$ be the set of vertices of $K$ that belong, respectively do not
  belong, to $K'$. We define a retraction on $K'$ as follows:
  \[\rho_s((x_0,g_1,\dots,x_n),(t_0,\dots,t_n))=((x_0,g_1,\dots,x_n),(u_0(s),\dots,u_n(s))),\]
  where
  \[u_i(s)=\begin{cases}
      \left(1+s\frac{\sum_{j\in B}t_j}{\sum_{j\in A}t_j}\right)t_i& \text{ if the $i$th vertex of $K$ belongs to $A$},\\
      (1-s)t_i& \text{ if the $i$th vertex of $K$ belongs to $B$}.
    \end{cases}
  \]
  This retraction is compatible with the face maps and degeneracy
  maps.

  We check that $\rho_s$ is the identity on simplices in
  $|\tile(\Pi)|$, because for these simplices we have $B=\emptyset$.
\end{proof}

Within the limit space, we see a collection of copies of the standard
tile, as follows: for $w\in A^*$, let $\tile_w$ be the full
subgroupoid on the objects $A^{-\N}w$. The same construction can be
performed for the geometric standard tile: given $w\in|L|_0$, set
\begin{equation}\label{eq:wbar}
  \overline w=\setsuch{y\in\operatorname{Star}(w)}{d(y,w)\le
    d(z,w)\forall z\in\operatorname{Star}(w)\cap|L|_0}.
\end{equation}S
For any $T\subset|L|_0$, set $\overline T=\bigcup_{w\in T}\overline
w$. Then for $w\in A^*$ we define the geometric tile
$\overline{\tile_w}\subset|L(\Pi)|$ as $\overline{A^{-\N}w}$.

\begin{prop}\label{prop:Ltiling}
  Assume $\Pi$ is nuclear; then for any $n\in\N$, the tiles
  $\setsuch{\overline{\tile_w}}{w\in A^n}$ cover $|L(\Pi)|$, and two
  tiles $\overline{\tile_w},\overline{\tile_{w'}}$ overlap if and only
  if $w,w'$ are connected in the Schreier graph of $\langle\Pi\rangle$
  on $A^n$.
  
  Furthermore each tile is the closure of its interior, and these
  tiles overlap with empty-interior intersection.
  
  The tiling $\setsuch{\overline{\tile_w}}{w\in A^{n+1}}$ is a
  refinement of the tiling $\setsuch{\overline{\tile_w}}{w\in A^n}$,
  with the tile $\overline{\tile_{xw}}$ being contained in the tile
  $\overline{\tile_w}$.
\end{prop}
\begin{proof}
  If two tiles $\overline{\tile_w},\overline{\tile_{w'}}$ overlap,
  then there is a morphism $(\alpha,\phi)$ of $L(\Pi)$ in the $n$th
  preimage of the critical locus, i.e.\ satisfying
  $\phi(-n)\neq\varepsilon$. The Schreier graph contains an edge,
  labeled $\phi(-n)$, from $w$ to $w'$.

  Conversely, let $g:w\to w'$ be such an edge. By
  Lemma~\ref{lem:tauepi}, there exists $(\alpha,\phi)\in L(\Pi)$ with
  $\phi(-n)=g$ and $\alpha([-n,-1])=w$. The tiles $\overline{\tile_w}$
  and $\overline{\tile_{w'}}$ overlap on the middle of that edge.

  The second statement is obvious, because the interior of a tile
  $\overline{\tile_w}$ is obtained by replacing $\le$ by $<$ in the
  definition of $\overline w$.
\end{proof}

Note that the tilings of the topological space $L(\Pi)_\top$ is much
more complicated~\cite{bartholdi-g-n:fractal}. The automaton $\Pi$
satisfies the \emph{open set condition} if for every $q\in Q$ there
exists $a\in A^n$ with $\tau(a,q)=\varepsilon$. This is equivalent to
asking that every $g\in\langle\Pi\rangle$ have a trivial state.

\begin{prop}[\cite{nekrashevych:ssg}*{Corollary~3.3.7}]
  If $\Pi$ satisfies the open set condition, then every tile in
  $L(\Pi)_\top$ is the closure of its interior, and for every $n\ge0$
  the tiles in the $n$th subdivision of $L(\Pi)_\top$ have disjoint
  interior.
  
  On the other hand, if $\Pi$ does not have the open set condition,
  then for every $n$ large enough one can find a tile in the $n$th
  subdivision of $L(\Pi)_\top$ which is covered by the other tiles in
  the same subdivision.
\end{prop}

Next, consider the solenoid $S(\Pi)$. This is naturally a foliated
space: for $\orbit$ a $\Gamma$-orbit on $A^\N$, let $F_\orbit$ denote
the full subgroupoid on the objects $A^{-\N}\orbit$. Note that there
are no morphisms in $S(\Pi)$ between distinct leaves.

We endow the leaf $F_\orbit$ with the ``left-Tychonoff,
right-discrete'' topology generated by open sets $\mathcal
O_{n,w}=\setsuch{\alpha:\Z\to A}{\alpha(n+i)=w_i\forall i\in\N}$ for
all $n\in\Z$ and $w\in A^\N$.

\begin{lem}
  Assume $\Pi$ is recurrent and spherically transitive. Then the shift
  map $f:S(\Pi)\to S(\Pi)$ preserves the foliation of $S(\Pi)$.
\end{lem}
\begin{proof}
  It suffices to show that the partition of $A^\Z$ in $A^{-\N}\orbit$
  is shift-invariant. For that purpose, we show that $\orbit$ is a
  $\Gamma$-orbit if and only if $A\orbit$ is a $\Gamma$-orbit.

  Pick $g\in\Gamma$, $a\in A$ and $v\in\orbit$. Then $g\cdot
  av=\sigma(a,g)\tau(a,g)\cdot v\in A\orbit$. Conversely, pick $a,b\in
  A$ and $v,w\in\orbit$. Then there exists $g\in\Gamma$ with $g\cdot
  a=b$, whence $g\cdot av=bv'$ with $v'\in\orbit$, and $g'\in\Gamma$
  with $g'\cdot v'=w$. Since $\Gamma$ is recurrent, there exists
  $h\in\Gamma_b$ with $\tau(b,h)=g'$; therefore $hg\cdot av=bw$.
\end{proof}

The leaves of the solenoid are again tiled spaces: for $w\in
\orbit\subset A^\N$, let $\tile_w$ be the full subgroupoid on the
objects $A^{-\N}w$. Define the geometric tiles $\overline{\tile_w}$ as
in~\eqref{eq:wbar}.
\begin{thm}
  Assume $\Pi$ is nuclear, and let $\orbit$ be a $\Gamma$-orbit in
  $A^\N$. Then $|F_\orbit|\subset |S(\Pi)|$ is tiled by
  $\setsuch{\overline{\tile_w}}{w\in\orbit}$.

  For any $w\in A^\N$, the adjacency graph of tiles on $|F_{\Gamma
  w}|$ is the Schreier graph of $\langle\Pi\rangle/\stab(w)$.

  Furthermore each tile is the closure of its interior, and these
  tiles overlap with empty-interior intersection.

  The shift map $f:S(\Pi)\to S(\Pi)$ sends each leaf $F_\orbit$ to
  another one carrying a refinement of the tiling of $F_\orbit$.  The
  tile $\overline{\tile_{xw}}\subset F_\orbit$ is mapped into the tile
  $\overline{\tile_w}$.
\end{thm}
\begin{proof}
  Analogous to that of Proposition~\ref{prop:Ltiling}.
\end{proof}

Every leaf $F_\orbit$ naturally covers $L(\Pi)$. The covering map
$F_\orbit\to L(\Pi)$ is given by the restriction of the natural
projection map $S(\Pi)\to L(\Pi)$ given by keeping only the negative
part of objects and morphisms. This covering map folds every tile
$\overline{\tile_w}$ of $|F_\orbit|$ onto $|L(\Pi)|$.

We note that $S(\Pi)\to L(\Pi)$ has the structure of a foliated bundle
$S(\Pi)=\tilde L\times_\Gamma A^\N$ in the sense
of~\cite{moore-c:foliated}. To see this, define the groupoid $\tilde
L$ as follows: its objects are $A^{-\N}\times\Gamma$.  Its set of
morphisms is the set of $(\alpha,\phi,g)\in
A^{-\N}\times\Gamma^{-\N}\times\Gamma$ such that $(\alpha,\phi)\in
L(\Pi)$ and $\phi_0=g$. There is a natural action of $\Gamma$ on
$A^\N$, and $\Gamma$ acts on $\tilde L$ by acting on the second
coordinate of $A^{-\N}\times\Gamma$.

In that context, we recall the notion of holonomy
groupoid~\cite{moore-c:foliated}*{Page~58}, restricted to a
transversal. Given a foliated space $S$ over $L$, fix a point $w\in
L$, and let $F$ be its fiber. The \emph{holonomy groupoid} has objects
$F$, and has a morphism from $f\in F$ to $f'\in F$ for every path
$\alpha$ in $S$ from $f$ to $f'$. Two paths $\alpha,\alpha'$ are
identified if they have the same holonomy, i.e.\ if both start at $f$,
end at $f'$, and $\alpha(\alpha')^{-1}$ induces by parallel transport
the identity of $F$ in the neighbourhood of $f$.

\begin{prop}
  The holonomy groupoid of the foliation, restricted to a transversal,
  is the groupoid of germs $\germ(\Pi)$ of $\Pi$.
\end{prop}
\begin{proof}
  Fix a point $w\in L(\Pi)$; then its fiber in $S(\Pi)$ is $wA^\N$.
  The parallel transport along $\alpha$ in a leaf of $S$ corresponds
  to the action of $\Gamma$ on $A^\N\cong wA^\N$. Two paths
  $\alpha,\alpha':wv\to wv'$ are equivalent if and only if
  $\alpha(\alpha')^{-1}$ is in the stable homotopy group of $v$.
\end{proof}

%Furthermore, there exists a maximal leaf, such that every leaf is a
%quotient of the standard leaf:
%\begin{prop}
%  Pick any $w\in L(\Pi)$, and let $\leaf$ be the cover of $L(\Pi)$
%  associated to the stable isotropy of $v$. Then the leaf associated
%  to the orbit of any $w\in A^\N$ is $\leaf/\germ(v)$.
%\end{prop}

%%%%%%%%%%%%%%%%%%%%%%%%%%%%%%%%%%%%%%%%%%%%%%%%%%%%%%%%%%%%%%%%
\section{Automata from orbispaces}\label{sec:a<o}
Let $L$ be a connected orbispace with a $d$-to-$1$ covering map
$f:L\to L$.

First, we may obtain a profinite group $\overline\Gamma$ as follows:
pick a base point $*\in L$. For any $n\in\N$, set
$A_n=f^{-n}(*)$. Consider the space $L_n=\{(g_0\in L,g_1\in
f^{-1}(g_0),\dots,g_n\in f^{-1}(g_{n-1}))\}$. Then $L\cong L_n$, by
projection on the last coordinate, and we may also view $L_n$ as a
degree-$d^n$ covering space over $L$, with $f^n:L_n\to L$ realized as
projection on the first coordinate.

Also, $f^n$ can be viewed as a bundle with structure group the
isometry group of an $n$-level $\#A$-regular rooted tree.  Consider
the associated principal bundle
\[Y_n=\bigg\{\Big(g\in L,p:\bigcup_{i=0}^nf^{-i}(g)\to\bigcup_{i=0}^n
A_i\text{ a rooted tree isometry}\Big)\bigg\},
\]
with covering map $\tilde f_n$ the projection on the first
coordinate; let $Z_n$ be the connected component of $(*,1,\dots,1)$ in
$Y_n$.  Then $Z_n$ is a Galois covering, with Galois group $\Gamma_n$.
Therefore $\tilde f:Y_{n+1}\to Y_n$ given by
$(g,p)\mapsto(g,p_{|\bigcup_{i=0}^nf^{-i}(g)})$ induces a projective
system $\tilde f:\Gamma_{n+1}\to\Gamma_n$. Set
$\overline\Gamma(L)=\varprojlim\Gamma_n$.

\begin{thm}\label{thm:sameclosure}
  Let $\Pi$ be nuclear and spherically transitive. Then
  $\overline\Gamma(L(\Pi))$ is the closure of $\langle\Pi\rangle$ in
  $\aut(A^*)$.
\end{thm}
\begin{proof}
  First, we remember by Lemma~\ref{lem:closure} that the closure of
  $\langle\Pi\rangle$ in $\aut(A^*)$ is the inverse limit of its
  finite quotients acting on $A^n$. This finite quotient is nothing
  but the permutation group $\Gamma_n$ in its action on $A_n$.
\end{proof}

We now seek a discrete group associated with $L$, and assume the
geometric realization $|L|$ of $L$ is arcwise connected. By
Proposition~\ref{prop:extend}, the map $f$ induces a $d$-to-$1$
covering map $|f|:|L|\to|L|$.

We place ourselves in the following situation, that of an arcwise
connected space $X$ endowed with a $d$-to-$1$ covering map $f:X\to X$.
Some of the interesting examples however, come from a slightly more
general situation, where $f$ needs only be a branched covering. This
reduces to the previous situation by removing from $X$ the branching
locus, as well as all its iterated direct and inverse images.

We actually do not need to remove the inverse images of the branching
locus; if we remove the forward images, we are led to consider a space
$X$ with a map $f$ defined on a dense subset of $X$, and satisfying
the unicity of path lifting property.

Here is the procedure for constructing an automaton $\Pi(X)$ out of
our data~\cite{nekrashevych:ssg}*{Proposition~5.2.2}. Pick a base point
$x\in X$ and choose for all $y\in f^{-1}(x)$ a path $\ell_y$ from $x$
to $y$. Let $K\triangleleft \pi_1(X,x)$ be the subgroup consisting of
all paths that induce the identity permutation on $f^{-n}(x)$, for all
$n\in\N$.  Our automaton is given as follows: its set of states $Q$ is
a subset of $\pi_1(X,x)$ that generates $\pi_1(X,x)/K$, and satisfies
the condition~\eqref{eq:cond} below. Its alphabet is $A=f^{-1}(x)$.
Given $a\in A$ and $q\in Q$, consider the preimage $\gamma$ of the
path $q$ with $\gamma(0)=a$.  The output function of $\Pi(X)$ is
$\sigma(a,q)=\gamma(1)$, and its transition function maps $(a,q)$ to a
path $\tau(a,q)\in Q$ congruent mod $K$ to
$\ell_a\gamma\ell_{\gamma(1)}^{-1}$. We therefore require for all
$a\in A$ and $q\in Q$:
\begin{equation}\label{eq:cond}
  \ell_a\gamma\ell_{\gamma(1)}^{-1}K\cap Q\neq\emptyset.
\end{equation}

Note that we may on one hand take $Q=\pi_1(X,x)$; however we wish
usually the set of states to be finite, in which case~\eqref{eq:cond}
expresses a non-trivial restriction on $Q$.

\begin{defn}[\cite{nekrashevych:ssg}*{Chapter~5}]\label{defn:img}
  We call $\Pi(X)$ the automaton constructed above, and, for $L$ a
  groupoid, we set $\Pi(L)=\Pi(|L|)$. The automaton $\Pi(X)$
  constructed depends of course on the choices made but we will show
  that the associated group $\langle\Pi(X)\rangle$ does not.
  
  The group $\langle\Pi(X)\rangle$ is called the \emph{iterated
    monodromy group} of $f$, written $\mathsf{IMG}(f)$.
\end{defn}

The following proposition shows that the automaton associated with a
\emph{topologically expanding} map is contracting. A metric version
was already proven in~\cite{nekrashevych:ssg}*{Theorem~5.5.3}.
\begin{defn}
  Let $X$ be a topological space and $f:X\to X$ a continuous map. It
  is \emph{topologically expanding} if there exists an open subset
  $\U\subset X\times X$ such that $\bigcap_{n\ge0}f^{-n}(\U)=\Delta$,
  where $\Delta=\setsuch{(x,x)}{x\in X}$ is the diagonal and $f$
  extends naturally to a function $X\times X\to X\times X$.

  If $Y\subset X\times X$ contains the diagonal, let us write $Y_0$
  for the connected component of $Y$ containing the diagonal. We then
  say $f$ is \emph{smoothly topologically expanding} if there exists a
  $\U\subset X\times X$ as above, satisfying furthermore
  $f^{-1}(\U)_0\subset\U$.
\end{defn}

\begin{prop}
  If $X$ is compact, locally arcwise connected, and $f$ is a smoothly
  topologically expanding degree-$d$ cover, then the automaton
  $\Pi(X)$ is contracting.
\end{prop}
\begin{proof}
  Set $A=\{1,\dots,d\}$.  Let $T$ be the standard $d$-regular rooted
  tree: it is the simplicial realization of the graph with vertices
  $A^*$ and an edge between $w$ and $wa$ for all $w\in A^*,a\in A$.
  The standard compactification of $T$ is $\overline T=T\cup\partial
  T$, with $\partial T=A^\N$.

  We assume the automaton $\Pi(X)$ has been constructed, with choices
  of basepoint $*$ and connecting paths $\{\ell_a\}$. The paths
  $\ell_a$ extend to paths $\ell_w$ starting at $*$, defined as
  follows. For $w\in A^*,a\in A$, set
  $\ell_{wa}=\ell_w\tilde\ell_{w,a}$, where $\tilde\ell_{w,a}$ is the
  $|w|$-th iterated $f$-preimage of $\ell_a$ starting at $\ell_w(1)$.
  This defines an embedding $T\to X$, mapping the vertex $w\in T$ to
  the extremity of $\ell_w$. Furthermore, for each $w\in A^\N$ we have
  a path $\ell_w:\R_{\ge0}\to X$ starting at $*$.

  We show that these paths are actually finite, i.e.\ the map $T\to X$
  extends to the compactification $\partial T\to X$. For this purpose,
  let $Y\subset X\times\partial T$ be the space of accumulation points
  of rays $\ell_w$:
  \[Y = \setsuch{(x,w)\in X\times\partial T}{x=\lim\ell_w(t_i)\text{
      for a sequence }t_i\to\infty}.\]
  Then $Y$ fibers over $\partial T$, with connected fibres since the
  $\ell_w$ are connected and $X$ is compact. We set
  \[Z = (Y\times_{\partial T}Y)\cap(\U\times\partial T).\]
  By assumption, $Z$ is stable under
  $f\times f\times\text{`shift'}$, so since $f$ is expanding
  $Z\subset\Delta\times\partial T$. It follows that the fibres of
  $Y\to\partial T$ are discrete, so $Y\cong\partial T$. We have shown
  that each ray $\ell_w$ has a unique accumulation point
  $\lim(\ell_w)$. This proves that the map $T\to X$ extends to
  $\overline T$.

  The relation $(x,y)\in\U$ means that ``$x$ and $y$ are close''. We
  need stronger notions of closeness, for which we introduce the
  notation $(x,y)\in\frac1k\U$, for $k\in\N$. By definition, if $\V$
  is an open neighbourhood of $\Delta$ in $X\times X$, then
  $\frac1k\V$ denotes a choice of an open subset of $X\times X$
  containing the diagonal and satisfying the condition:
  \[\text{ if }(x_i,x_{i+1})\in\frac1k\U\text{ whenever }0\le
  i<k,\text{ then }(x_0,x_k)\in\V.
  \]
  Let us show that such sets always exist. For this, fix $\V$ and $k$,
  and let $(Y_n)$ be a decreasing sequence (or net) of compact
  neighbourhoods of the diagonal $\Delta$ satisfying $\bigcap
  Y_n=\Delta$. Then $p_{12}^{-1}(Y_n)\cap \dots\cap
  p_{k,k+1}^{-1}(Y_n)\searrow \Delta^{(k+1)}$, where $\Delta^{(k+1)}$
  is the diagonal in $X^{k+1}$, and $p_{ij}$ are the projections on
  two coordinates.  Since images commute with decreasing intersections
  of compact sets, we get $Z_n=p_{1,k+1}(p_{12}^{-1}(Y_n)\cap\dots\cap
  p_{k,k+1}^{-1}(Y_n)) \searrow \Delta$. For $n$ big enough, we will
  then have $Z_n\subset \V$, and we set $\frac1k\V = (Y_n)^\circ$ for
  such a choice of $n$.

  If $\gamma:[0,L]\to X$ is a path and $\V$ is an open neighbourhood
  of $\Delta$, we will write $\gamma\Subset\V$ to mean
  $(\gamma(t),\gamma(t'))\in\V$ for all $t,t'\in[0,L]$; more
  generally, if $Y\subset X$, we write $Y\Subset\V$ to mean
  $(y,y')\in\V$ for all $y,y'\in Y$.

  If $w\in T$, the \emph{cone} $C(w)$ of $w$ is the image in $X$ of
  the subset spanned by $wA^*\cup wA^\N$ in $\overline T$.

  Set $\V=\frac13\frac12\U$.  We may choose $R\in\N$ such that
  $C(w)\Subset\V$ for all $w\in A^R$. To construct such an $R$,
  consider the function $h:\partial T\to\N$,
  \[h(w)\triangleq\min\setsuch{M\in\N}{C(w_1\dots w_M)\Subset\V}.\]
  This function is continuous on a compact, so is bounded. Set $R=\max
  h(\partial T)$.

  Define also $N,N',K'\subset\pi_1(X,*)$ by
  \begin{align*}
    N &= \setsuch{\ell_v\rho\ell_w^{-1}}{|v|,|w|\ge R\text{ and
      }\rho\Subset\V},\\
    N'&= \setsuch{\ell_v\rho\ell_w^{-1}}{|v|=|w|=R\text{ and
      }\rho\Subset\frac12\U},\\
    K'&= \setsuch{\ell_w\rho\ell_w^{-1}}{\rho\Subset\U}.
  \end{align*}

  First, we claim that the elements of $K'$ act trivially on the tree
  of preimages of $*$, and therefore are trivial in
  $\langle\Pi(X)\rangle$.  Take $\gamma\in K'$. Then $\gamma$ is a
  ``balloon'' $\rho\Subset\U$, attached by a ``string'' $\ell_w$ to
  $*$. Since $f^{-1}(\U)_0\subset\U$, an $f$-preimage of $\rho$ is
  again of the form ``balloon in $\U$ attached by a string''. They are
  all loops, so $\gamma$ acts trivially by monodromy on $f^{-1}(*)$.
  Iterating this procedure, we see that $\gamma$ acts trivially on
  $f^{-n}(*)$ for all $n\in\N$.

  We next claim that $N\subset N'$. Take $\gamma\in N$,
  $\gamma=\ell_v\rho\ell_w^{-1}$. Set $v'=v_1\dots v_R$ and
  $w'=w_1\dots w_R$. Let $\ell_{v',v}$ be the image in $X$ of the path
  from $v'$ to $v$ in $T$, and similarly for $\ell_{w',w}$. Then
  $\ell_{v',v},\rho,\ell_{w',w}^{-1}\Subset\V$, so
  $\rho'\triangleq\ell_{v',v}\rho\ell_{w',w}^{-1}\Subset\frac12\U$, and we
  can write $\gamma=\ell_{v'}\rho'\ell_{w'}^{-1}\in N'$.

  We claim that $N'$ is finite modulo $K'$. Indeed take two elements
  $\gamma,\gamma'\in N'$ with the same $v,w$:
  $\gamma=\ell_v\rho\ell_w^{-1},\gamma'=\ell_v\rho'\ell_w^{-1}$. Their
  quotient is a balloon $\ell_v\rho(\rho')^{-1}\ell_v^{-1}$. It is in
  $K'$ since $\rho(\rho')^{-1}\Subset\U$. We are left with finitely
  many choices for $v,w\in A^R$.

  Finally, we claim that for any $\gamma\in\pi_1(X,*)$ there exists
  $n\ge R$ such that all $f^{-n}$-preimages of $\gamma$ are contained
  in $\V$. Indeed partition $\gamma$ in small segments
  $\gamma_1,\dots,\gamma_k$ such that $\gamma_i\Subset\U$ for all
  $i\in\{1,\dots,k\}$. Since $f^{-n}(\U)$ converges to $\Delta$, there
  exists $n\in\N$ such that $f^{-n}(\U)\Subset\frac1k\V$. The preimages
  of $\gamma_i$ are in $\frac1k\V$, so all preimages of $\gamma$ belong
  to $\V$, and therefore define an element of $N$.

  We have shown that $N$ contains the nucleus and that its image in
  $\langle\Pi(X)\rangle$ is finite.
\end{proof}
Note that this proposition applies in particular if $f\in\C(z)$ and
the postcritical orbit of $f$ does not intersect the Julia set of
$f$; one then takes for $X$ the Julia set of $f$~\cite{fatou:funceq}.

Let $A,A'$ be two sets of same cardinality. In the following
definition, we relax the definition of automata in allowing the input
and output alphabets to be respectively $A$ and $A'$. More precisely,
the transition function remains a function $\tau:A\times Q\to Q$, but
the output function becomes a function $\sigma:A\times Q\to A'$. Such
an automaton, with an initial state $q_0\in Q$, defines a tree
homomorphism $A^*\to(A')^*$, which is an isomorphism precisely when
$\sigma(-,q)$ is a bijection for all $q\in Q$. The
\emph{states}\footnote{By abuse of notation, we refer to the states of
  $\phi$ to mean the states of an automaton defining the map $\phi$.}
of $\phi:A^*\to(A')^*$ are the maps $\phi':A^*\to(A')^*$ given by
$\phi(vw)=v'\phi'(w)$ for some $v\in A^*$ and any $w\in A^*$. The set
of states of the map $q_0:A^*\to(A')^*$ is a subset of $Q$.

\begin{defn}\label{defn:equivalent}
  Two automata $\Pi$ and $\Pi'$ on alphabets $A$ and $A'$ are
  \emph{equivalent} if $\#A=\#A'$, and there exists a tree isomorphism
  $\phi:A^*\to (A')^*$ given by a finite automaton in the sense
  described above, such that the states of $\phi$ are all of the form
  $\phi g$ for some $g\in\langle\Pi\rangle$, and
  \[\phi\circ\langle\Pi\rangle\circ\phi^{-1}=\langle\Pi'\rangle.\]
\end{defn}

\begin{prop}\label{prop:treeiso}
  Let $(X,f)$ be a space with a covering and let $\Pi$, $\Pi'$ be two
  automata constructed as above from possibly different data
  $x,\ell_y,Q$ and $x',\ell'_y,Q'$. Assume $\Pi$ is contracting. Then
  $\Pi$ and $\Pi'$ are equivalent.
\end{prop}
\begin{proof}
  Since equivalence of automata is transitive, we may assume
  $Q\subset\pi_1(X,x)$ and $Q'=\pi_1(X,x')$. Choose a path $s$ from
  $x$ to $x'$ and use it to identify $\pi_1(X,x)$ with $\pi_1(X,x')$.
  We get an injection
  \[s_Q:\langle Q\rangle\subset\pi_1(X,x)\stackrel{\sim}{\longrightarrow} 
    \pi_1(X,x')=\langle Q'\rangle.
  \]
  The path $s$ also induces a bijection between the preimages
  $A=f^{-1}(x)$ and $A'=f^{-1}(x')$, and similarly between
  higher-order preimages of $x$ and $x'$, by lifting appropriately $s$
  through powers of $f$. The paths $\ell_y$ can be used to identify
  the tree of preimages of $x$ with $A^*$ and similarly to identify
  the preimages of $x'$ with $(A')^*$. Composing these three tree
  isomorphisms gives an isomorphism $s_A:A^*\to (A')^*$.

  For all $a\in A$, let $s_a$ be the $f$-preimage of $s$ starting at
  $a$, and call its other extremity $a'\in A'$. Define
  $g_a\in\pi_1(X,x)$ by
  \[g_a=\ell_as_a\ell_{a'}^{-1}s^{-1}.\] Define recursively
  $\phi:A^*\to (A')^*$ by $\phi(aw)=a'\phi(g_aw)$.  Then $\phi=s_A$,
  and since $\Pi$ is contracting, $\phi$ is automatic
  by~\cite{nekrashevych:ssg}*{Corollary~2.11.7}.  Furthermore, the
  states of $\phi$ are clearly of the form $\phi g$ for some
  $g\in\langle\Pi\rangle$; for instance, on the first level, they are
  precisely the $\{\phi g_a\}$.

  The best way to understand the definition of $\phi$ is as follows:
  forget for an instant that $s$ is fixed, and denote the resulting
  $\phi$ by $\phi_s$. Then
  \[\phi_s\circ g_a=\phi_{g_as}=\phi_{\ell_a s_a\ell_{a'}^{-1}},\]
  in accordance the transition function defined above~\eqref{eq:cond}.

  Clearly the actions of $\langle Q\rangle$ and $\langle Q'\rangle$ on
  $A^*$ and $(A')^*$ are intertwined by $s_Q$ and $s_A$. We are left to
  show that the induced map
  $s_\Pi:\langle\Pi\rangle\to\langle\Pi'\rangle$ is an isomorphism,
  where $\langle\Pi\rangle$ and $\langle\Pi'\rangle$ are the images of
  $\langle Q\rangle$ and $\langle Q'\rangle$ in $\aut(A^*)$ and
  $\aut((A')^*)$ respectively.  The injectivity of $s_\Pi$ follows from
  that of $s_Q$.  To show the surjectivity of $s_\Pi$, let
  $g'\in\aut((A')^*)$ be in $\langle\Pi'\rangle$, let $\gamma'\in\langle
  Q'\rangle=\pi_1(X,x')$ be a path representing it, and set
  $\gamma=s\gamma's^{-1}\in\pi_1(X,x)$. We need to find a word in $Q$
  that acts on $A^*$ the same way that $\gamma$ does; such a word
  exists because $Q$ was assumed to generate the image of $\pi_1(X,x)$
  in $\aut(A^*)$.
\end{proof}

Consider a manifold $X$ with a branched covering map $f$, and let $P$
be the smallest closed subset of $X$ containing the critical values of
$f$ and $f(P)$. If $P$ does not disconnect $X$ we may consider the
space $X'=X\setminus P$, with a partially defined covering map $\tilde
f:X'\dashrightarrow X'$; we still have the unique lifting property for
loops via $\tilde f$, so the construction of $\mathsf{IMG}(f)$ is
unaffected by the fact that $\tilde f$ is not defined everywhere.

We may also consider the space
$X''=X\setminus\overline{\bigcup_{m,n\in\N}f^{-m}f^n(\text{critical
    points})}$, now with an everywhere defined covering map. An
important source of examples is given by $X$ the Riemann sphere, and
$f$ a rational map --- see for instance~\cite{bartholdi-g-n:fractal}.

\begin{thm}\label{thm:samegp}
  If $\Pi$ is nuclear and smooth, then $\Pi$ is equivalent to
  $\Pi(L(\Pi))$; more precisely, the data $x,\ell_y,Q$ may be chosen
  in Definition~\ref{defn:img} so that $\Pi=\Pi(L(\Pi))$.
\end{thm}

\begin{lem}\label{lem:tverse>word}
  Pick $x\in|L(\Pi)|_0$ and $\gamma\in\pi_1(|L(\Pi)|_1,x)$. Assume
  that $\gamma$ is transverse to the midpoints of the critical edges.
  Let $u\in Q^*$ be the sequence of last states $\phi_0$ of the
  critical edges $(\alpha,\phi)$ crossed by $\gamma$. Then the action
  of $\gamma$ on the tree of preimages of $x$ only depends on $u$.
\end{lem}
\begin{proof}
  Let $n>0$ be an integer. We want to show that the action of $\gamma$
  on $f^{-n}(x)$ only depends on $u$. Let $y\in f^{-n}(x)$ be a
  preimage and $\tilde\gamma$ be the unique lift of $\gamma$ starting
  at $y$. As before, $|L(\Pi)|$ has a tiling by tiles
  $\overline{\tile_w}$, for all $w\in A^n$. The points where
  $\tilde\gamma$ crosses from a tile to another are exactly the
  preimages of the points where $\gamma$ crosses a critical edge.  So
  if $t\in[0,1]$ is one of these points and
  $\tilde\gamma(t-\epsilon)\in\overline{\tile_w}$ for small
  $\epsilon>0$, then
  $\tilde\gamma(t+\epsilon)\in\overline{\tile_{\sigma(w,q)}}$, where
  $q=\phi_0$ is the last state of the critical edge $(\alpha,\phi)$ on
  which $\gamma(t)$ lies.  By induction, we get
  $\tilde\gamma(1)\in\overline{\tile_{\sigma(w,u)}}$ where $w\in A^n$
  is such that $y=\tilde\gamma(0)$ belongs to $\overline{\tile_w}$.
  This shows that $\tilde\gamma(1)$ is the unique preimage of $x$
  lying in $\overline{\tile_{\sigma(w,u)}}$ and is therefore
  determined by $u$ and $w$ only.
\end{proof}

\begin{proof}[Proof of Theorem~\ref{thm:samegp}]
  We show that a judicious choice of base point, connecting paths and
  generating set for $|L(\Pi)|$ gives $\Pi(|L(\Pi)|)\simeq\Pi$.

  We start with any base point $x\in|L(\Pi)|_0$. For every $q\in
  Q\setminus\{1\}$, let $(\alpha_q,\phi_q):\alpha_q\to\beta_q$ be a
  critical edge with $(\phi_q)_0=q$. Such an edge can be constructed
  inductively using the expansion rule: $(\alpha_q)_{n-1}\triangleq
  e_{(\phi_q)_n}$, $(\phi_q)_{n-1}\triangleq v_{(\phi_q)_n}$. Let
  $\gamma^1_q:x\to\alpha_q$ and $\gamma^2_q:\beta_q\to x$ be paths
  that cross no critical edge.  These paths are images of paths in the
  standard tile, which is connected by
  Proposition~\ref{prop:tileconn}. We let $\tilde Q$ be the set of
  paths
  $\gamma_q\triangleq\gamma^1_q(\alpha_q,\phi_q)\gamma^2_q\subset\pi_1(|L(\Pi)|,x)$,
  parametrized by $q\in Q$. The $\ell_y$'s are taken to be paths from
  $x$ to $y\in f^{-1}(x)$ that don't cross any critical edge. These
  paths give the natural identification of the tree of preimages of
  $x$ with $A^*$, namely $y\in f^{-n}(x)\mapsto y_{-n}y_{-n+1}\ldots
  y_{-1}\in A^*$.
  
  We show that the $\gamma_q$ generate the image of
  $\pi_1(|L(\Pi)|,x)$ in $\aut(A^*)$. Let $\gamma$ be any loop. By
  transversality, we may assume that it satisfies the hypothesis of
  Lemma~\ref{lem:tverse>word}. We get a word $w=w_1\ldots w_r\in Q^*$
  and by Lemma~\ref{lem:tverse>word}, $\gamma_{w_1}\ldots\gamma_{w_r}$
  acts the same way as $\gamma$ on $A^*$, which is what we wanted.
  
  So we only need to show that $\Pi(|L(\Pi)|)\simeq(A,Q)$ has same
  output and transition functions as $\Pi$. We claim that if $y$ is
  the preimage of $x$ corresponding to $a$, then $\sigma(y,\gamma_q)$
  is the preimage of $x$ corresponding to $\sigma(a,q)$. Letting
  $\tilde\gamma$ be the lift of $\gamma_q$ starting at $y$, the
  argument of Lemma~\ref{lem:tverse>word} shows
  $\sigma(y,\gamma_q)=\tilde\gamma(1)\in\overline{\tile_{\sigma(a,q)}}$,
  which proves our claim.
  
  Last, we need to show that
  $\ell_y\tilde\gamma\ell_{\tilde\gamma(1)}^{-1}$ acts the same way on
  $A^*$ as $\tau(a,q)$ does. By Lemma~\ref{lem:tverse>word} it is
  enough to show the following:
  $\ell_y\tilde\gamma\ell_{\tilde\gamma(1)}^{-1}$ has at most one
  critical edge. If it has one, its last state is $\tau(a,q)$. And if
  it doesn't then $\tau(a,q)=1$. Recall that the $\ell_y$'s were
  chosen without critical edges and that $\tilde\gamma$ is a preimage
  of $\gamma$, which has a unique critical edge $(\alpha_q,\phi_q)$.
  So the only possible critical edge of
  $\ell_y\tilde\gamma\ell_{\tilde\gamma(1)}^{-1}$ is the preimage of
  $(\alpha_q,\phi_q)$. Let us call $(\tilde\alpha,\tilde\phi)$ that
  preimage. We need to show $\tilde\phi_0=\tau(a,q)$; this is clear
  since $\tilde\phi_{-1}=(\phi_q)_0=q$ and $\tilde\alpha$ is in the
  same tile $\overline{\tile_a}$ as $y$, which implies
  $\tilde\alpha_{-1}=a$.

  If $\Pi(|L(\Pi)|)$ had been built using different data, then by
  Proposition~\ref{prop:treeiso}, we would still have $\Pi$ equivalent
  to $\Pi(|L(\Pi)|)$.
\end{proof}

\begin{thm}
  If $\Pi$ and $\Pi'$ are nuclear and equivalent, then $(L(\Pi),f)$
  and $(L(\Pi'),f')$ are Morita equivalent.
  
  If furthermore $\Pi$ and $\Pi'$ are smooth, then $\Pi$ and $\Pi'$
  are equivalent if and only if $(L(\Pi),f)$ and $(L(\Pi'),f')$ are
  Morita equivalent.
\end{thm}
\begin{proof}
  Let $\phi_0:A^*\to(A')^*$ be an equivalence between $\Pi$ and
  $\Pi'$. Denote $\langle\Pi\rangle$ by $\Gamma$, and set
  $\Phi=\phi_0\Gamma$. The Morita equivalence is given by
  \begin{multline*}
    P=\left\{(\alpha,\phi,\alpha')\in
      A^{-\N}\times\Phi^{-\N}\times(A')^{-\N}\right|
    \phi(n)(\alpha(n)w)=\alpha'(n)\phi(n+1)(w)\\
    \text{ for all }w\in A^\N\text{ and }n<-1,
    \text{ and }\phi(-\N)\text{ is finite}\big\},
  \end{multline*}
  with source and target maps $s_P(\alpha,\phi,\alpha')=\alpha$ and
  $t_P(\alpha,\phi,\alpha')=\alpha'$.
  
  We first show that $P$ is compact; this is because
  $W=\bigcup_{(\alpha,\phi,\alpha')\in P}\phi(-\N)$ is finite; we
  actually show that $W$ is contained in the Cartesian product of the
  set of states of $\phi$ and the nucleus of $\Pi$. Fix some
  $p=(\alpha,\phi,\alpha')\in P$, and consider the set of $g\in\Gamma$
  that satisfy $\phi_0g=\phi(n)$ for some $n<0$. Since $\Pi$ is
  nuclear, for every such $g$ there exists a $k\in\N$ with
  $\tau(A^k,g)\subseteq Q$. Since there are finitely many such $g$'s
  associated with $p$, there is a common such $k$ for all $g$'s.  Fix
  now $n<0$, and write $\phi(n)=\tau(w,\phi(n-k))$ with
  $w=\alpha(n-k)\dots\alpha(n-1)$. Write $\phi(n-k)=\phi_0g'$ for a
  $g'\in\Gamma$. Then
  $\phi(n)=\tau(w,\phi_0g')=\tau(\sigma(w,g'),\phi_0)\tau(w,g')$. We
  observe that $\tau(\sigma(w,g'),\phi_0)$ is a state of $\phi_0$, and
  that $\tau(w,g')\in\tau(A^k,g')\subseteq Q$, which proves that $W$
  is finite; therefore $P$ is closed in the compact $(A\times W\times
  A')^{-\N}$, so is compact.
  
  We next show that $s_P$ and $t_P$ are onto. Choose any $\alpha\in
  A^{-\N}$.  Then for every $m<0$, there exists
  $(\alpha,\phi,\alpha')\in A^{-\N}\times\Phi^{-\N}\times(A')^{-\N}$
  satisfying $\phi(n)(\alpha(n)w)=\alpha'(n)\phi(n+1)(w)$ for all
  $w\in A^\N$ and $n\ge m$; for instance, set $\phi(m)=\phi_0$. This
  determines $\phi(n)$ and $\alpha'(n)$ for all $n\ge m$. Choose the
  $\phi(n),\alpha'(n)$ arbitrarily for $n<m$. Since $P$ is compact,
  there exists an accumulation point of the above choices as
  $m\to-\infty$, and $s_P$ is onto. The same argument applies to
  $t_P$.
  
  We then show that $L(\Pi)\sqcup P\sqcup P^{-1}\sqcup L(\Pi')$ is a
  groupoid. It suffices to check that given $x\in L(\Pi)$ and $y,z\in
  P$ with $t(x)=s(y)$, the product $xy\in P$ is well-defined; if
  $t(y)=t(z)$ then $yz^{-1}\in L(\Pi)$ is well-defined; and if
  $s(y)=s(z)$ then $y^{-1}z\in L(\Pi')$ is well-defined. Write
  $x=(\alpha,\phi)$ and $y=(\beta,\psi,\beta')$. Then
  $xy=(\alpha,\phi\psi,\beta')$ and $(\phi\psi)(-\N)$ is finite
  because both $\phi(-\N)$ and $\psi(-\N)$ are finite, so $xy\in P$.
  
  Similarly, write $z=(\alpha,\phi,\beta')$. Then
  $yz^{-1}=(\beta,\psi\phi^{-1},\alpha)$; and
  $\psi(n)\phi(n)^{-1}=\phi_0g_1g_2^{-1}\phi_0^{-1}\in\Gamma$ assumes
  finitely many values, so $yz^{-1}\in L(\Pi)$.
  
  Write also $z=(\beta,\phi,\alpha')$. Then
  $y^{-1}z=(\beta',\psi^{-1}\phi,\alpha')$; and
  $\psi(n)^{-1}\phi(n)=g_1^{-1}g_2\in\Gamma$ assumes finitely many
  values, so $y^{-1}z\in L(\Pi)$.
  
  Finally, the shift map clearly extends to $P$, by deleting the
  $(-1)$st entry of $\alpha$, $\phi$ and $\alpha'$.
  
  Conversely, assume that $\Pi$ and $\Pi'$ are nuclear and smooth, and
  let $P$ be a Morita equivalence between $G=L(\Pi)$ and $G'=L(\Pi')$.
  Set $G''=G\sqcup P\sqcup P^{-1}\sqcup G'$. We obtain a diagram
  $|G|\hookrightarrow|G''|\hookleftarrow|G'|$, and the inclusions
  commute with the coverings $f$, $f'$ and $f''$. By the second part
  of Theorem~\ref{thm:samegp}, we may pick data $x,\ell_y,Q$ such
  that $\Pi(|G|)=\Pi$, and similarly data $x',\ell_{y'},Q'$ such that
  $\Pi(|G'|)=\Pi'$. We may then push these data forward into $|G''|$,
  obtaining two sets of data again giving $\Pi$ and $\Pi'$
  respectively. It now follows from Proposition~\ref{prop:treeiso}
  that $\Pi$ and $\Pi'$ are equivalent automata.
\end{proof}

\section{Examples}
This section describes some examples of automata and their associated
limit spaces. We start by exhibiting various automata that fail to
satisfy the various conditions: contraction, smoothness, etc. We then
describe the favourable situation of the automaton associated xxxx to the
covering $f(z)=z^2-1$ of the Riemann sphere.

We will describe the group $G$ of the automaton by giving its
decomposition map $\psi:G\to G\wr\sym A$ on generators of $G$. The
associated automaton can be recovered by taking as states the
generators of $G$, as alphabet $A$, and taking a transition from state
$q$ to state $q'$ with input $a$ and output $a'$ precisely when
$\psi(q)=(r,\pi)$ with $r(a)=q'$ and $\pi(a)=a'$.

The automaton correponding to this description can be obtained by
drawing square tiles with labels the input and output alphabet letters
and states, by the procedure described in Subsection~\ref{subs:automata}.

\subsection{The Lamplighter group}
Here, as in the next two examples, $A=\{0,1\}$ and $\sigma=(0,1)$ is
the non-trivial element of $\sym A$.

The Lamplighter group is the group
$G=(\Z/2)\wr\Z=\oplus_\Z(\Z/2)\rtimes\Z$. It may be generated by an
automaton as follows: $G=\langle a,b\rangle$ with $\psi(a)=(a,b)$ and
$\psi(b)=(b,a)\sigma$. This automaton is not nuclear; indeed it is not
even contracting, since $a$ has infinite order and the projection of
$\psi^n(a)$ on the first vertex is $a$. Its square description is as
follows:
\[\begin{matrix}
  \sqtile 0a0a & \sqtile 1a1b & \sqtile 0b1b & \sqtile 1b0a
\end{matrix}
\]

This automaton representation of $G$ was used by Grigorchuk and \.Zuk
in~\cite{grigorchuk-z:lamplighter} to compute the $\ell^2$ spectrum of
the simple random walk on $G$.

The identification of $G$ with the automaton can be understood as
follows: identify the boundary of the tree $A^\N$ with
$\F[2][[t]]$, under the map $(a_i)\mapsto\sum a_it^i$. Then $a$ and
$b$ identify respectively with the affine maps $f\mapsto (1+t)f$ and
$f\mapsto (1+t)f+1$ of $\F[2][[t]]$. Therefore $G$ identifies with the
maps of the form $f\mapsto (1+t)^nf+p$ for some
$n\in\Z$ and $p\in\F[2][1+t,(1+t)^{-1}]$.

\subsection{The Baumslag-Solitar group}
Let $m,n$ be two integers. The Baumslag-Solitar group $G_{m,n}$ is
defined by its presentation
\[G_{m,n}=\langle a,t|t^{-1}a^mt=a^n\rangle.\]
It is residually finite precisely when $m=\pm n$ or $m=\pm1$ or
$n=\pm1$.
If $m=1$, it may be represented by affine transformations as $a(X)=X+1$
and $t(X)=nX$.

The group $G_{1,3}$ may be generated by an automaton as follows:
$G=\langle a,b,c\rangle$ with $\psi(a)=(a,b)$ and
$\psi(b)=(a,c)\sigma$ and $\psi(c)=(b,c)$. Again this automaton is not
nuclear. Its square description is as follows:
\[\sqtile 0a0a\quad \sqtile 1a1b\quad \sqtile 0b1a\quad \sqtile 1b0c\quad \sqtile 0c0b\quad \sqtile 1c1c\]

The identification of $G$ with the automaton can be understood as
follows: identify the boundary of the tree $A^\N$ with $\Z_2$,
under the map $(a_i)\mapsto\sum a_i2^i$. Then $a$, $b$ and $c$
identify respectively with the affine maps $X\mapsto 3X$, $X\mapsto
3X+1$ and $X\mapsto 3X+2$. Therefore $G$ identifies with the maps of
the form $X\mapsto 3^nf+p$ for some $n\in\Z$ and $p\in\Z[1/3]$.

The similarity between the lamplighter and Baumslag-Solitar groups is
not accidental; a common construction of $G_{1,m}$ and $(\Z/q)\wr\Z$
by automata is described in~\cite{bartholdi-s:solvable}.

\subsection{The Odometer}
Again identify the boundary $A^\N$ of the tree with $\Z_2$, and
consider the subgroup $\Z$ of $\Z_2$. This cyclic group may be
generated by an automaton as follows: $G=\langle\tau\rangle$ with
$\psi(\tau)=(\varepsilon,\tau)\sigma$, where as before $\varepsilon$
denotes the identity state. The associated automaton is contracting,
with nucleus $\{\varepsilon,\tau,\tau^{-1}\}$. Its square description
is as follows:
\[\begin{matrix}
  \sqtile 0\tau1\varepsilon & \sqtile 1\tau0\tau & \sqtile 0{\tau^{-1}}1{\tau^{-1}} & \sqtile 1{\tau^{-1}}0\varepsilon
\end{matrix}
\]

The limit space $L(\Pi)$ is $I/\{0^\N=1^\N\}$ with $I$ as in
Lemma~\ref{lem:interval}, and its topological quotient $L(\Pi)_\top$
is homeomorphic to the circle.  The standard tile is $I$ and its
topological quotient is $[0,1]$. The topological quotient of the
associated solenoid is the standard $2$-adic solenoid: the inverse
limit of
\[\xymatrix{{\cdots}\ar[r] & {S^1}\ar[r]^{()^2} & {S^1}\ar[r]^{()^2} & {S^1}}.\]

\vspace*{1cm}
We consider in the next three examples some contracting actions which
exhibit various ``pathologies''.
\subsection{A non-recurrent example}
Take now $A=\{0,1,2\}$, with $\sigma=(0,1,2)$ a three-cycle, and
consider the action of $\Z$ defined as follows:
$G=\langle\tau\rangle$, with
$\psi(\tau)=(\varepsilon,\tau,\tau)\sigma$.  The associated automaton
is contracting, with nucleus
$\{\varepsilon,\tau^{\pm1},\tau^{\pm2}\}$.  The square description of
its nucleus is as follows (we omit the squares for the state
$\varepsilon$, which have vertical labels $\varepsilon$ and equal
horizontal labels):
\begin{gather*}
  \begin{matrix}
    \sqtile 0{\tau^{-2}}1{\tau^{-2}} &
    \sqtile 0{\tau^2}2\tau &
    \sqtile 0{\tau^{-1}}2{\tau^{-1}} &
    \sqtile 0\tau1\varepsilon \\
    \sqtile 1{\tau^{-2}}2{\tau^{-1}} &
    \sqtile 1{\tau^2}0{\tau^2} &
    \sqtile 1{\tau^{-1}}0\varepsilon &
    \sqtile 1\tau2\tau \\
    \sqtile 2{\tau^{-2}}0{\tau^{-1}} &
    \sqtile 2{\tau^2}1\tau &
    \sqtile 2{\tau^{-1}}1{\tau^{-1}} &
    \sqtile 2\tau0\tau
  \end{matrix}
\end{gather*}

This group is not recurrent: indeed the stabilizer of a vertex (say
$0$) is $\langle\tau^3\rangle$, and its projection on the subtree
$0A^*$ is $\langle\tau^2\rangle$.

The limit space $L(\Pi)_\top$ is the standard $2$-adic solenoid. In
particular, it is connected, but not arcwise connected. Its
self-covering is the ``triple-the-angle'' map.

\subsection{A non-smooth example}
Take again $A=\{0,1,2\}$ and $\sigma=(0,1,2)$, and consider the action
of $\Z$ defined as follows: $G=\langle\tau\rangle$, with
$\psi(\tau)=(\tau,\tau^{-1},\tau)\sigma$. The associated automaton is
contracting, with nucleus $\{\varepsilon,\tau^{\pm1},\tau^{\pm2}\}$.  The square
description of its nucleus is as follows:
\begin{gather*}
  \begin{matrix}
    \sqtile 0{\tau^{-2}}1\varepsilon &
    \sqtile 0{\tau^2}2\varepsilon &
    \sqtile 0{\tau^{-1}}2{\tau^{-1}} &
    \sqtile 0\tau1\tau \\
    \sqtile 1{\tau^{-2}}2{\tau^{-2}} &
    \sqtile 1{\tau^2}0\varepsilon &
    \sqtile 1{\tau^{-1}}0{\tau^{-1}} &
    \sqtile 1\tau2{\tau^{-1}} \\
    \sqtile 2{\tau^{-2}}0\varepsilon &
    \sqtile 2{\tau^2}1{\tau^2} &
    \sqtile 2{\tau^{-1}}1{\tau} &
    \sqtile 2\tau0\tau
  \end{matrix}
\end{gather*}

The stabilizer of a vertex (say $0$) is $\langle\tau^3\rangle$, which
projects to $G$ on the subtree $0A^*$; so $G$ is recurrent.

The minimal automaton generating $G$, with set of states
$Q=\{\varepsilon,\tau,\tau^{-1}\}$, is not smooth. However, since
$\psi(\tau^2)=(\varepsilon,\varepsilon,\tau^2)\sigma^2$, the subgroup
$\langle\tau^2\rangle$ is smooth.

The limit space $L(\Pi)_\top$ is a circle, but $|L(\Pi)|$ looks more
like a M\"obius strip. The three tiles of $|L(\Pi)|$ project to the
three overlapping intervals $[0,2\pi/3]$, $[2\pi/3,\pi/3]$ and
$[\pi/3,0]$ of the circle:
\begin{center}
  \includegraphics{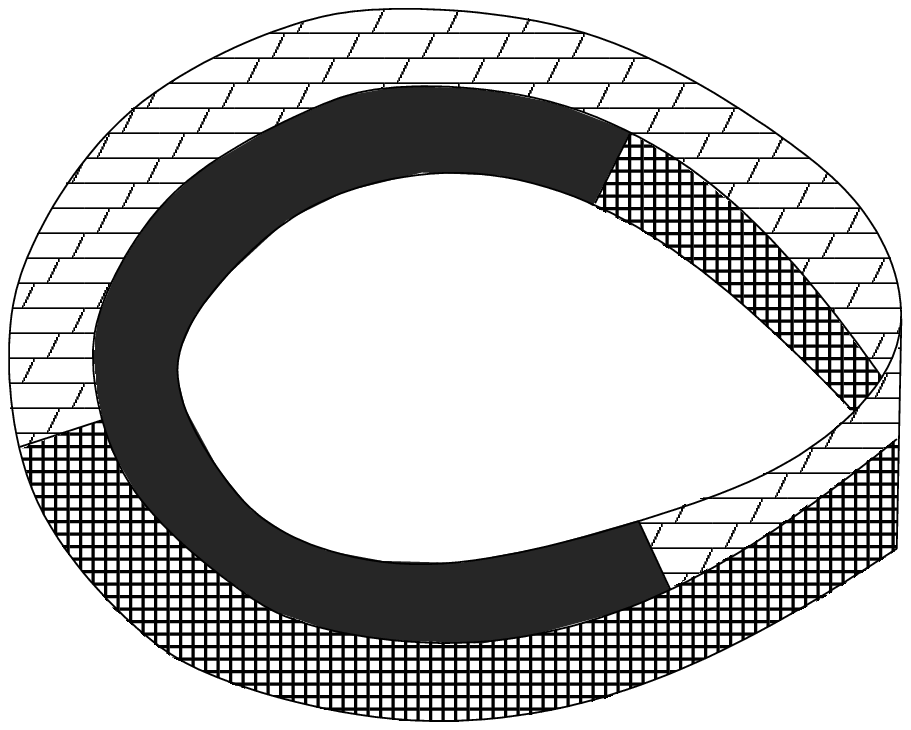}
\end{center}

\subsection{A more complicated non-smooth example}
Take again $A=\{0,1,2\}$ and $\sigma=(0,1,2)$, and consider the action
of $\Z$ defined as follows: $G=\langle\tau\rangle$, with
$\psi(\tau)=(\tau^2,\varepsilon,\tau^{-1})\sigma$. The associated automaton is
contracting, with nucleus $\{\varepsilon,\tau^{\pm1},\tau^{\pm2}\}$.  The square
description of its nucleus is as follows:
\begin{gather*}
  \begin{matrix}
    \sqtile 0{\tau^{-2}}1\tau &
    \sqtile 0{\tau^2}2{\tau^2} &
    \sqtile 0{\tau^{-1}}2{\tau} &
    \sqtile 0\tau1{\tau^2} \\
    \sqtile 1{\tau^{-2}}2{\tau^{-1}} &
    \sqtile 1{\tau^2}0{\tau^{-1}} &
    \sqtile 1{\tau^{-1}}0{\tau^{-2}} &
    \sqtile 1\tau2\varepsilon \\
    \sqtile 2{\tau^{-2}}0{\tau^{-2}} &
    \sqtile 2{\tau^2}1{\tau} &
    \sqtile 2{\tau^{-1}}1\varepsilon &
    \sqtile 2\tau0{\tau^{-1}}
  \end{matrix}
\end{gather*}

The stabilizer of a vertex (say $0$) is $\langle\tau^3\rangle$, which
projects to $G$ on the subtree $0A^*$; so $G$ is recurrent.

The automaton generating $G$ is again non-smooth, but this time in an
essential way: there is no element of the nucleus sending $0$ to $1$
and projecting to a power of the trivial state. Indeed the only
element with that property is $\tau^{-5}$ and, since it does not
belong to the nucleus, there is no (bounded-length) path in $|L(\Pi)|$
between $w0$ and $w'1$ for any $w,w'\in A^{-\N}$.

The limit space $L(\Pi)$ is therefore connected, but not arcwise
connected. Its topological quotient is as before a circle, since the
group is recurrent. This example illustrates how the property of being
arcwise connected is not invariant under Morita equivalence.

The topological quotient of the tile is the closure of the union of
countably many closed intervals. Its boundary is a Cantor set.

Below we represent the tiling of the circle (drawn as a horizontal
interval) by three copies of the standard tile. We cross each copy of
the tile by an interval in order to make the picture more legible, and
draw each copy in a different tint.
\begin{center}
  \includegraphics{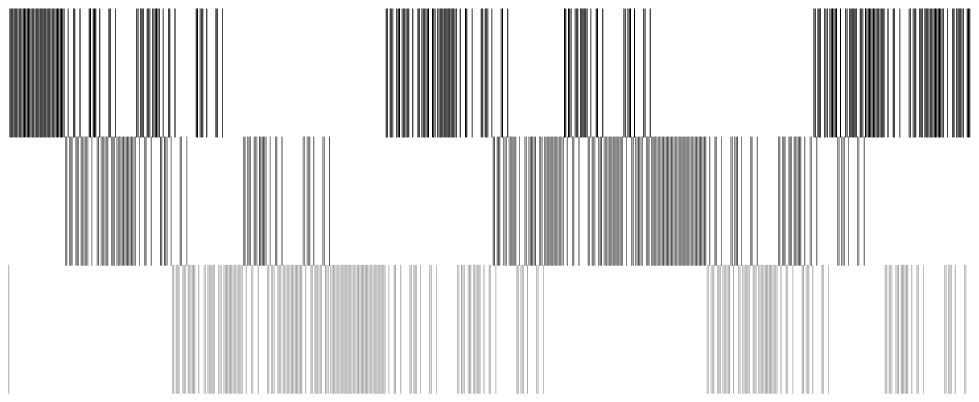}
\end{center}

\subsection{The ``basilica group''}
This highly non-trivial example of group served as a motivation for
the study of iterated monodromy groups and their general properties.
It is defined, with $A=\{0,1\}$ and $\sigma=(0,1)$, as follows:
$G=\langle a,b\rangle$, with $\psi(a)=(\varepsilon,b)\sigma$ and
$\psi(b)=(\varepsilon,a)$. The associated automaton is contracting, with nucleus
$\{\varepsilon,a^{\pm1},b^{\pm1},a^{-1}b,b^{-1}a\}$.

The group $G$ is torsion-free, amenable, but cannot be obtained from
groups of subexponential growth via direct limits, extensions,
subgroups and quotients~\cite{bartholdi-v:amenability}. It was the first
example with such a property. More details on $G$ can be found
in~\cite{grigorchuk-z:torsionfree}.

The topological quotient of the limit space is the Julia set of the
complex map $f(z)=z^2-1$.

The odometer is conjugate, within $\aut A^*$, to the subgroup $\langle
a^{-1}b\rangle$ of $G$. This explains that the limit space of $G$ is a
quotient of the limit space of the odometer.

\end{document}